\documentclass[12pt,reqno,a4paper]{amsart}

\usepackage[a4paper, left=3cm, right=2.5cm, top=2.5cm, bottom=2.5cm]{geometry}

\usepackage{standalone}

\usepackage[backend=bibtex, style=alphabetic, sorting=nyt]{biblatex}
\renewbibmacro{in:}{}

\usepackage{graphicx}
\usepackage{xcolor}
\usepackage{tikz}
\usepackage{placeins}
\usepackage{float}
\usetikzlibrary{arrows.meta, positioning}
\graphicspath{{figures/}}

\usepackage{amsmath, amssymb}
\usepackage{amsthm}
\usepackage{textcomp}
\usepackage[OT1]{fontenc}
\usepackage{stmaryrd}
\usepackage{enumitem}
\setlist[enumerate]{label= (\arabic*) }

\usepackage{booktabs}
\usepackage{array}

\usepackage[colorlinks=true, linkcolor=black, citecolor=black, urlcolor=blue]{hyperref}

\usepackage{caption,subcaption}

\usepackage{listings}
\usepackage{extarrows}
\usetikzlibrary{braids,knots, hobby, spath3, decorations.markings, arrows.meta, positioning, decorations.pathreplacing}

\theoremstyle{plain}
\newtheorem{theorem}{Theorem}[section]
\newtheorem{lemma}[theorem]{Lemma}
\newtheorem{corollary}[theorem]{Corollary}

\newtheorem*{theorem*}{Theorem}

\theoremstyle{remark}
\newtheorem*{remark}{Remark}

\theoremstyle{definition}
\newtheorem{definition}[theorem]{Definition}
\newtheorem{example}[theorem]{Example}

\newcommand{\sqbinom}[2]{\genfrac{[}{]}{0pt}{}{#1}{#2}}

\makeatletter
\renewcommand{\@tocwrite}[2]{%
  \ifx\@secnumber\@empty
  \else
    \@xp\@tocwriteb\csname toc#1\endcsname{#1}{#2}%
  \fi
}
\makeatother

\title{Colored HOMFLY homology for positive two-strand torus knots}
\author{Yuan Xu}
\date{\today}

\begin{document}
\begin{abstract}
    The triply graded colored HOMFLY homology, as a categorification of the colored HOMFLY polynomial, is a powerful colored link invariant. In this paper, we study the reduced column-colored HOMFLY homology of a family of positive two-strand torus knots. For every pair of integers $m\geq 1$ and $p\geq 0$, we give an explicit formula for the bottom row, i.e.\ the $A=0$ specialization, of the reduced column-colored HOMFLY homology $\overline{\mathrm{HHH}}_{\bigwedge^m}(T(2,2p+1))$ of the positive torus knot $T(2,2p+1)$. This formula verifies a conjecture of Gorsky, Gukov and Sto{\v{s}}i{\'c}~\cite{gorsky2013quadruply}.
    
    The proof is computational. Following the recursive formulae for the reduced column-colored HOMFLY homology of positive torus links due to Conners~\cite{conners2025doubly}, the graded dimension of bottom row of the reduced column-colored homology of $T(2,2p+1)$ can be expressed as a formal power series, and simplified through an elimination process encoded by $\kappa$-admissible sequences to a summation. The summands consist of exponents of $T$, $Q$ and $1-Q$, where the exponents of $1-Q$ are determined by an orbit tracing argument, while the exponents of $Q$ and $T$ are computed by analyzing the numbers and the positions of ones in the $\kappa$-admissible sequences. Applying a quantum multinomial theorem, the powers of $Q$ can be assembled into a product of quantum binomials depending only on the number of ones in the admissible sequences, which yields the explicit formula and verifies the conjecture in~\cite{gorsky2013quadruply}.
\end{abstract}    
\maketitle
\tableofcontents
\section{Introduction}

\subsection{Knot polynomials and their categorification}
Ever since the invention of the Jones polynomial in 1984~\cite{jones1985polynomial}, knot polynomials and their categorification have been a main tool to distinguish knots and links. Building on the discovery of various knot polynomials including the Alexander polynomial~\cite{alexander1928topological}, $\mathfrak{sl}_N$ polynomial~\cite{reshetikhin1991invariants}, and HOMFLY polynomial~\cite{freyd1985new}, Khovanov~\cite{khovanov2000categorification} introduces a bigraded link homology theory using matrix factorization, whose Euler characteristic degenerates to Jones polynomial. This homology theory is now called \emph{Khovanov homology}, and such an extension is called a ``\emph{categorification}''.

Later categorifications of other knot polynomials were introduced. For instance, knot Floer homology categorifies the Alexander polynomial~\cite{rasmussen2003floer,ozsvath2004holomorphic}, and Khovanov-Rozansky homology categorifies $\mathfrak{sl}_N$ polynomials~\cite{khovanov2004matrix}, which recovers Khovanov homology when $N=2$. Unlike the one-variable polynomials above, the two-variable HOMFLY polynomial is categorified to a triply graded homology called HOMFLY homology, first introduced by Khovanov and Rozansky~\cite{khovanov2008matrix} and later reformulated by Khovanov~\cite{khovanov2007triply}.

\subsection{Colored invariants}
All the knot invariants above can be understood, in the language of representation theory, as assigning every link component a standard vector representation $V$ of $U_q(\mathfrak{sl}_N)$. In this sense, a natural generalization is to assign an arbitrary irreducible representation $V^{\lambda}$, indexed by a \emph{Young diagram} $\lambda$, to every link component. Such a choice of $\lambda$ for every link component is called a \emph{coloring}, and invariants stemming from these constructions are called \emph{colored invariants}. 

The rigorous construction of such triply-graded colored homology evolves over years after the discovery of the (uncolored) HOMFLY homology. Mackaay, Sto{\v{s}}i{\'c}, and Vaz~\cite{mackaay20111} first defined a $(1,2)$-colored HOMFLY homology, and proved that it is a link invariant, while for general $\bigwedge^k$-colored homology they only gave a potential definition without a proof of invariance. The existence and invariance for \emph{column-colored homology} was later proved using geometric methods by Webster and Williamson~\cite{webster2017geometric}. The $\lambda$-colored HOMFLY homology was constructed by Cautis~\cite{cautis2017remarks}, for arbitrary Young diagram $\lambda$ and arbitrary links. In particular, when $\lambda=(1^k)$ is a column with a single 1 at each row, the resulting colored homology will be called \emph{the column-colored HOMFLY homology}. 

In practice, we will work with \emph{the reduced colored HOMFLY homology} of a knot, which is finite-dimensional  and recovers all the information that the unreduced colored homology encodes. 

\subsection{On computation of triply graded link homology}
Although the existence and invariance of colored HOMFLY homology theory has been well established, the direct computation of link homology still remains challenging. 
For the uncolored version, Elias, Hogancamp, and Mellit~\cite{elias2019computation,hogancamp2019torus,mellit2022homology} gave a recursive formula for the uncolored HOMFLY homology of positive torus links $T(m,n)$, using the inductive construction of the categorical Young symmetrizer $P_n$. 
The computational result in uncolored HOMFLY homology is later extended to the \emph{reduced column-colored} HOMFLY homology of positive torus links $T(m,n)$ with coloring $(\bigwedge^k,\bigwedge^1,\ldots,\bigwedge^1)$ by Conners~\cite{conners2025doubly}, where a family of recursive formulae is given. This paper deals exclusively with reduced column-colored HOMFLY homology. 

Although the recursive formulae for colored torus link homology are now available, reducing them to an explicit formula remains a nontrivial problem. Gorsky, Gukov, and Sto{\v{s}}i{\'c}~\cite{gorsky2013quadruply} conjectured an explicit formula for the bottom row of the reduced column-colored HOMFLY homology of the positive torus knot $T(2,2p+1)$, based on an analogous identification established for other version of link homologies.
This paper mainly focuses on a family of torus knots $T(2,2p+1)$ with $\bigwedge^m$-coloring, and studies the $A=0$ specialization (called \emph{bottom row} for simplicity) of the reduced HOMFLY homology, which is given as the following theorem.
\begin{theorem*}[Main result]
    Let $m\geq 1$, $p\geq 0$ be two integers. Then the bottom row of reduced column-colored HOMFLY homology of $\bigwedge^m$-colored positive torus knot $T(2,2p+1)$ has dimension 
        \[
            \dim \left( \overline{\mathrm{HHH}}_{\bigwedge^m} (T(2,2p+1)) \right) |_{A=0} = 
            Q^{\frac{3m^2-3m}{4}}
            \mathop{\sum\sum\cdots\sum}_{0\leq k_p\leq k_{p-1}\leq\cdots\leq k_1\leq k_0=m}
            \sqbinom{m}{k_1}_{Q^{-1}}\sqbinom{k_1}{k_2}_{Q^{-1}}\times\cdots
        \]
        \[
            \times\sqbinom{k_{p-1}}{k_p}_{Q^{-1}} Q^{\sum\limits_{j=1}^{p}k_{j-1}k_j-2m\sum\limits_{j=1}^{p}k_j}T^{\sum\limits_{j=1}^{p}k_j},
        \]
        where $A, Q, T$ are the three gradings of the reduced column-colored HOMFLY homology.
\end{theorem*}
This result verifies the conjecture in~\cite{gorsky2013quadruply}.

\subsection{Organization of the paper}
The rest of this paper is organized as follows. 

\textbf{Section~\ref{chapter:algebraic-setup}} reviews necessary notions to define the reduced column-colored homology. 
At the beginning, we introduce the 2-category $\mathrm{SSBim}$ of singular Soergel bimodules based on the building blocks \emph{merges} ${}_{a+b}M_{a,b}$ and \emph{splits} ${}_{a,b}S_{a+b}$. 
Then, we define the \emph{Rickard complex} $C_{\mathrm{KR}}(\beta)$ associated to colored braid $\beta$, and the \emph{reduced Rickard complex $C_{\mathrm{KR}}(\beta,p)$ at a marked point $p$}, and then review the Hochschild homology functor $\mathrm{HH}_*$ for singular Soergel bimodules and the reduced Hochschild homology functor $\mathrm{HH}_{\neg p}$ for complexes of singular Soergel bimodules. 
In the end of this section, we prove a version of \emph{quantum multinomial equation}, which is a generalization of the multinomial equation to the quantum integers. This equation is useful in expressing the final result. 

\textbf{Section~\ref{chapter:computation}} contains the core computations. To begin with, using theorems in~\cite{conners2025doubly}, the dimension of the column-colored homology of the torus knot $T(m,n)$ can be written as 
\[
    \dim\left(\overline{\mathrm{HHH}}_{\bigwedge^k}\right) (T(m,n)) = (1-T) {f(k)}^{-1} p_e(v_m, w_m)
\]
where $f(k)$ is a quantum constant determined by $k$, tuple $(v_m,w_m) := (1^m 0^m, 1^m 0^{2pm})$ (to be defined precisely in Section~\ref{chapter:computation}), permutation $e=e_l\in\mathfrak{S}_l$ is the identity permutation in $l$ variables, and $p_\sigma(v,w)$ is a formal power series in $(A,Q,T)$ satisfying the recursive formula, with permutation $\sigma\in\mathfrak{S}_l$, and $|v|=|w|=l$.
To describe the convention to simplify the power series $p_e$, we introduce the notion of $\kappa$-admissible sequences. By applying the recursive formulae with a certain pattern, the formal power series $p_e(v_k,w_k)$ can be expressed as a sum over all $\kappa$-admissible sequence $(u_0,u_1,u_2,\ldots,u_p)$, where the summand only consists of $Q$, $1-Q$, and $T$. 
By an orbit tracing theorem in~\cite{conners2025doubly}, the exponent of $(1-Q)$ is exactly the number of orbit changes of a permutation ${C(v,w)}_\sigma$ associated to the triple $(\sigma,v,w)$, which only depends on the number of ones in the sequence. By direct computation, the exponent of $Q$ and $T$ can be traced by analyzing the numbers and the positions of ones in the $v$-entry (or $w$-entry), where the exponent of $T$ depends also only on the number of ones. Applying the quantum multinomial equation in Section~\ref{chapter:algebraic-setup}, the summand in $Q$ can be written as a product of binomials (or a single multinomial). Eventually, this gives an explicit formula for the bottom row of the homology, which verifies the conjecture in~\cite{gorsky2013quadruply}.

\subsection{Conventions}
Here are some notations and conventions that will be used throughout this paper.
\begin{itemize}
    \item $\mathfrak{S}_l$: the permutation group in $l$ variables;
    \item ${[k]}_Q:=\frac{1-Q^k}{1-Q}$: the $k$-th quantum integer, where $k$ is a non-negative integer.
    \item ${[n]}_Q!:=\prod_{k=1}^{n}{[k]}_Q$: the $n$-th quantum factorial, where $n$ is a non-negative integer.
    \item $\sqbinom{n}{k}_Q:=\frac{{[n]}_Q!}{{[k]}_Q!{[n-k]}_Q!}$: the quantum binomial.
    \item If the variable $Q$ is clear from the context, we will omit the subscript $Q$ in the quantum integers, quantum factorials, and quantum binomials.
    \item $\mathbf{1}_A: U\to \{0,1\}$: the characteristic function of $A$, defined as 
    \[
        \mathbf{1}_A(x):=
        \begin{cases}
            1, &x\in A\\
            0, &x\notin A
        \end{cases}.
    \]
    \item For positive integers $m$ and $n$, we write 
    \[
        (1^m 0^n):=(\underbrace{11\ldots 1}_{m}\underbrace{00\ldots 0}_{n}),
        \]
        to denote the $(01)$-sequence of length $m+n$ with $m$ ones followed by $n$ zeros.
    \item Unless otherwise stated, the tensor products $\otimes$ are taken over the ground field $\mathbb{C}$, all diagrams of braids, webs, and permutations are oriented from right to left, and the strands are numbered from top to bottom. 
    \item Throughout this paper, we will use several independent grading groups. To distinguish among them, we assign a formal variable to each $\mathbb{Z}$-grading group, and write the degrees as products of powers of these formal variables. For example, if $x$ is an element of $q$-degree $2$ and $t$-degree $3$, then we write $\deg(x)=q^2 t^3$. A $\mathbb{Z}$-graded $R$-module $M$ with its formal variable $q$ is also called a \emph{$\mathbb{Z}_q$-graded $R$-module}. 
\end{itemize}
\subsection*{AI disclosure}
The author used ChatGPT 6 Astra and DeepSeek Flash V4.1 to verify simple cases of Lemma~\ref{section:quantum-multinomial-formula}, to generate a preliminary framework for Figure~\ref{fig:elimination-process}, and to polish the language. No AI tools were used for any other aspect of this work. The author takes full responsibility for the mathematical arguments and final text.
\subsection*{Acknowledgments}
This paper is also my master's thesis from the University of Zurich. I would like to thank my supervisor Anna Beliakova, for her overall suggestions on this paper and excellent lectures throughout my master's studies. Her guidance during my master's studies set me on the path towards research in this area of mathematics. 
Most of all, I would like to thank Luke Conners, who suggested the topic of this paper. Without his patient explanations of the theory and his careful review of my drafts, this paper would not have been possible. His guidance throughout the project and his constant encouragement have been invaluable.
Finally, I would like to thank my friends Yidan Tang and Kai Zhu for their support and encouragement throughout my studies. I would also like to thank my parents for their unwavering love and support, without which I would not have been able to pursue my studies and research in mathematics.

\section{Algebraic setup}\label{chapter:algebraic-setup}
In this section, we review the basic concepts for colored HOMFLY homology theory, and prove a quantum multinomial identity. Unless otherwise noted, strands in webs and braids are oriented from right to left, and the order of the strands are numbered from top to bottom. For example, the uppermost strand on the right is the first strand, and the lowermost strand on the right is the last strand. 

\subsection{Singular Soergel bimodules}\label{section:singular-soergel-bimodule}
    In this section, we follow~\cite[Section 3]{hogancamp2024link} to introduce the 2-category of singular Soergel bimodules in this section. To begin with, we introduce the symmetric polynomial algebra which acts on the singular Soergel bimodules.
    
    \begin{definition}[Symmetric polynomial algebra]
        For a non-negative integer $N$, define $R :=\mathbb{Q}[x_1,x_2,\ldots, x_N]$ as the $\mathbb{Z}_q$-graded ring with $\deg_q(x_i)=2$, for every $1\leq i\leq N$. If $\mathbf{a}=(a_1,a_2,\ldots,a_m)\vdash N$ is a composition of $N$, then 
        \[
            \mathfrak{S}_{\mathbf{a}}:= \mathfrak{S}_{a_1}\times \mathfrak{S}_{a_2}\times \cdots \times \mathfrak{S}_{a_m}
        \]
        acts on $R$ by permutation on the variables. Denote $R^{\mathbf{a}}=R^{(a_1,a_2,\ldots,a_m)}$ as the invariant subring of $R$ under the action of $\mathfrak{S}_{\mathbf{a}}$. If $\mathbb{X}_j = \{x_{a_{j-1}+1},\ldots, x_{a_j}\}$ are disjoint alphabets for $1\leq j\leq m$, then we also write 
        \[
            \mathrm{Sym}(\mathbb{X}_1|\mathbb{X}_2|\cdots|\mathbb{X}_m):=R^{\mathbf{a}}.
        \]
    \end{definition}

    Bimodules over symmetric polynomial algebras can be organized into a 2-category, which is defined as follows:

    \begin{definition}
        Let $N\geq 0$ be a non-negative integer. Define $\mathrm{Bim}_N$ as the 2-category with the following data: 
        \begin{itemize}
            \item objects are integer compositions $\mathbf{a}=(a_1,\ldots,a_m)$ of $N$, where $a_i\geq 1$, $1\leq i\leq m$;
            \item 1-morphisms from $\mathbf{a}$ to $\mathbf{b}$ are $(R^{\mathbf{b}},R^{\mathbf{a}})$-bimodules;
            \item 2-morphisms are graded bimodule homomorphisms,
        \end{itemize}
        where the composition $\star$ of 1-morphisms is given by the tensor product over the common part $R^{\mathbf{a}}$, and composition of 2-morphisms is given by composition of bimodule homomorphisms. Moreover, write $\mathrm{Bim}:=\sqcup_{N\geq 0} \mathrm{Bim}_N$.
    \end{definition}

    There is a monoidal structure on the 2-category $\mathrm{Bim}$ defined as follows.

    \begin{definition}[External tensor product]
        For positive integers $N_1$ and $N_2$, define 
        \[
            \boxtimes: \mathrm{Bim}_{N_1}\times \mathrm{Bim}_{N_2} \to \mathrm{Bim}_{N_1+N_2}
        \]
        by $\mathbf{a}\boxtimes\mathbf{b}:= (a_1,\ldots,a_m,b_1,\ldots,b_n)$ for objects $\mathbf{a}=(a_1,\ldots,a_m)\vdash N_1$ and $\mathbf{b}=(b_1,\ldots,b_n)\vdash N_2$, 
        \[
            M \boxtimes M' := M\otimes M', 
        \]
        and
        \[
            f\boxtimes g := f\otimes g,
        \]
        where $M$ and $M'$ are $(R^{\mathbf{b}}, R^{\mathbf{a}})$-bimodule and $(R^{\mathbf{b}'}, R^{\mathbf{a}'})$-bimodule, respectively, and $f$ and $g$ are bimodule homomorphisms.
    \end{definition}
    Note that the $\mathbb{Z}_q$-grading on $R$ induces a canonical $\mathbb{Z}_q$-grading on every $R$-bimodule. Besides the operations on bimodules, we introduce the following notation for graded bimodules.   
    \begin{definition}[Grading shifts]
        Suppose that $M$ is a $\mathbb{Z}_q$-graded bimodule and $n\in\mathbb{Z}$ is an integer. Write $M=\bigoplus_{d\in\mathbb{Z}}M_d$ as the homogeneous decomposition of $M$. Define $q^n M$ as the $\mathbb{Z}_q$-graded module with its degree $d$ component 
        \[
            {(q^n M)}_d:=M_{d-n},\quad \forall d\in \mathbb{Z}.
        \]
        Moreover, for a chain complex $C$ of $\mathbb{Z}_q$-graded bimodules, let $t$ denote its homological grading. Define $t^n C$ as the complex with its homological degree $d$ component
        \[
            {(t^n C)}_d := C_{d-n},\quad \forall d\in \mathbb{Z}.
        \]
    \end{definition}
    \begin{remark}
        Note that every graded bimodule $M$ can be viewed as a complex concentrated in homological degree $0$. In this sense, the grading shift $t^n$ can be also applied to graded bimodules to indicate its homological degree. For example, the complex 
        \[
            C=\left(
                0\to C_0 \to C_1 \to C_2\to 0
            \right)
        \]
        can be rewritten as 
        \[
            C=\left(
                0\to t^0 C_0 \to t^1 C_1 \to t^2 C_2\to 0
            \right).
        \]
    \end{remark}
    Singular Soergel bimodules are 1-morphisms generated under grading shifts, tensor products and external tensor products, by elementary building blocks called \emph{merge} and \emph{split}, which are defined as follows.
    \begin{definition}[Merge, split]
        Let $a$ and $b$ be two positive integers. Define the merge of $a$ and $b$ as the following $(R^{a+b},R^{(a,b)})$-bimodule
        \[
            {}_{a+b}M_{a,b}:=q^{-ab} R^{a+b}\otimes_{R^{a+b}} R^{(a,b)},
        \]
        and the split of $a$ and $b$ as the following $(R^{(a,b)},R^{a+b})$-bimodule 
        \[
            {}_{a,b}S_{a+b}:= R^{(a,b)}\otimes_{R^{a+b}} R^{a+b},
        \]
        where $R=\mathbb{C}[x_1,x_2,\ldots,x_{a+b}]$ is the polynomial ring of $(a+b)$ variables. Moreover, the merge and split can be graphically represented as certain trivalent graphs called \emph{webs} as demonstrated in Figure~\ref{fig:merge-and-split-in-numbers}.
        \begin{figure}[h]
            \centering
            \begin{subfigure}{0.45\textwidth}
                \centering
                \begin{tikzpicture}[
    baseline = (current bounding box.center),
    scale = 0.5,
    ultra thick,
]
    \draw (4,1) ..controls +(180:0.5) and +(90:0.5).. (3,0);
    \draw (3,0) ..controls +(-90:0.5) and +(180:0.5).. (4,-1);
    \draw (1,0) -- (3,0);

    \node[above] at (1,0) {${a+b}$};
    \node[right] at (4,1) {${a}$};
    \node[right] at (4,-1) {${b}$};
\end{tikzpicture}
                \caption{Merge of $a$ and $b$}
            \end{subfigure}
            \begin{subfigure}{0.45\textwidth}
                \centering
                \begin{tikzpicture}[
    baseline = (current bounding box.center),
    scale = 0.5,
    ultra thick,
]
    \draw (0,1) ..controls +(0:0.5) and +(90:0.5).. (1,0);
    \draw (1,0) ..controls +(-90:0.5) and +(0:0.5).. (0,-1);
    \draw (3,0) -- (1,0);

    \node[left] at (0,1) {${a}$};
    \node[left] at (0,-1) {${b}$};
    \node[above] at (3,0) {${a+b}$};
\end{tikzpicture}
                \caption{Split of $a$ and $b$}
            \end{subfigure}
            \caption{Merge and split}\label{fig:merge-and-split-in-numbers}
        \end{figure}
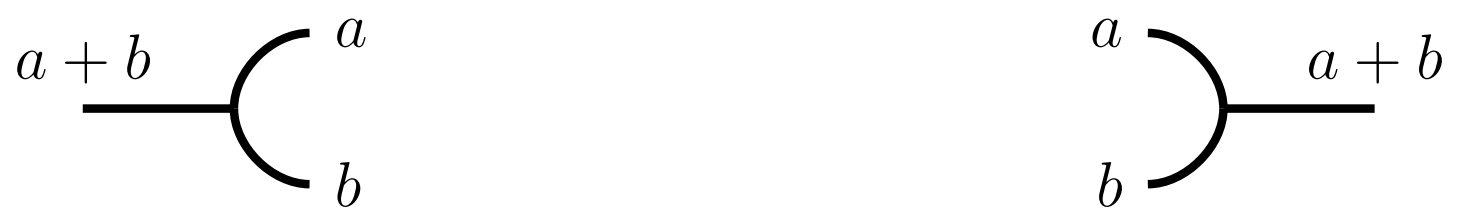
    \end{definition}
     
    \begin{remark}
        The size of an alphabet $\mathbb{X}$ is also denoted as $|\mathbb{X}|$. We often replace the labels on braids and webs by alphabets of variables. To be precise, suppose that $\mathbb{U}=\{u_1,\ldots,u_a\}$, $\mathbb{V}=\{v_1,\ldots,v_b\}$, and $\mathbb{W}=\{w_1,\ldots,w_{a+b}\}$ are disjoint alphabets, then we can replace the label $a$ and $b$ on the merge and split by $\mathbb{U}$, $\mathbb{V}$, and $\mathbb{W}$, respectively. 
        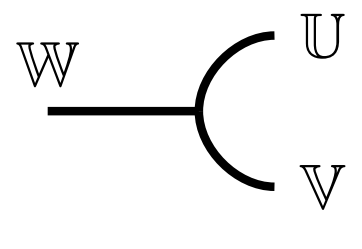
\begin{figure}[h]
            \centering
            \begin{tikzpicture}[
    baseline = (current bounding box.center),
    scale = 0.5,
    ultra thick,
]
    \draw (4,1) ..controls +(180:0.5) and +(90:0.5).. (3,0);
    \draw (3,0) ..controls +(-90:0.5) and +(180:0.5).. (4,-1);
    \draw (1,0) -- (3,0);

    \node[above] at (1,0) {$\mathbb{W}$};
    \node[right] at (4,1) {${\mathbb{U}}$};
    \node[right] at (4,-1) {${\mathbb{V}}$};
\end{tikzpicture}
            \caption{Merge in terms of alphabets $\mathbb{U}$, $\mathbb{V}$, and $\mathbb{W}$}\label{fig:merge-and-split-in-alphabets}
        \end{figure}
        \FloatBarrier{}
    \end{remark}

    Now we can define the category of singular Soergel bimodules as the full subcategory of $\mathrm{Bim}_N$ generated by merge and split. 
    \begin{definition}[Singular Soergel bimodules]\label{def:singular-soergel-bimodule}
        Let $N\geq 0$ be a non-negative integer. The 2-category $\mathrm{SSBim}$ is defined as the full subcategory of $\mathrm{Bim}$, with the same objects as $\mathrm{Bim}$, and 1-morphisms in $\mathrm{SSBim}$ are generated by ${}_{a+b}M_{a,b}$ and ${}_{a,b}S_{a+b}$ under grading shifts, taking direct summands, compositions $\star$, and external tensor products $\boxtimes$. 1-morphisms in $\mathrm{SSBim}$ are called \emph{singular Soergel bimodules}.
    \end{definition}

    Operations $\star$ and $\boxtimes$ are interpreted as horizontal composition and vertical stacking, respectively. Therefore, every singular Soergel bimodule has a graphical presentation. In the following context, we will use this graphical presentation to denote the singular Soergel bimodule associated to a web. 

    In addition to the 2-strand merge and split, the following example shows how vertical stacking plays a role in the construction of singular Soergel bimodules.
    \begin{example}
        Here is an example of a singular Soergel bimodule with more ends.
        \begin{align*}
            \vcenter
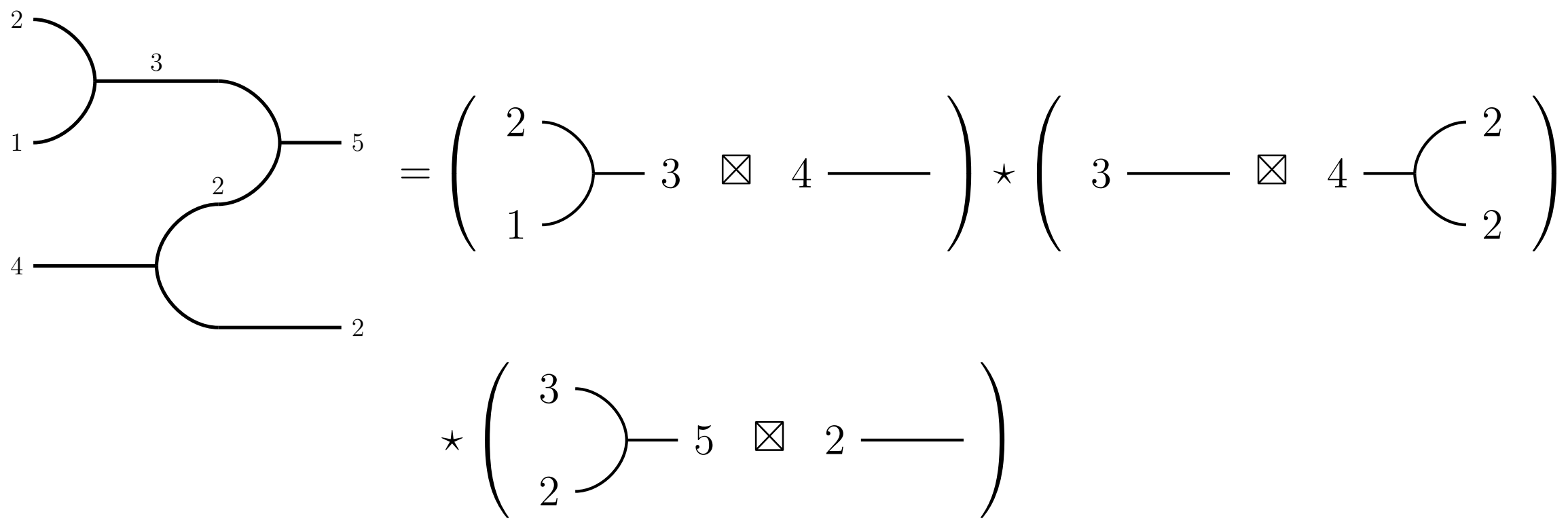
{\hbox{
                \includegraphics[scale = 0.6]{figures/webs/vertical-stacking-demo.tex}
            }}
                =&\left(\vcenter{\hbox{
                    \begin{tikzpicture}[
    baseline = (current bounding box.center),
    scale = 0.5,
    thick,
]
    \draw (0,1) ..controls +(0:0.5) and +(90:0.5).. (1,0);
    \draw (1,0) ..controls +(-90:0.5) and +(0:0.5).. (0,-1);
    \draw (2,0) -- (1,0);

    \node[left] at (0,1) {$2$};
    \node[left] at (0,-1) {$1$};
    \node[right] at (2,0) {$3$};
\end{tikzpicture}
                }}
                \boxtimes
                \vcenter{\hbox{
                    \begin{tikzpicture}[
                        baseline = (current bounding box.center),
                        thick,
                    ]
                        \draw (0,0) -- (1,0);
                        \node[left] at (0,0) {$4$};
                    \end{tikzpicture}
                }}
                \right)\star \left(
                    \vcenter{\hbox{
                    \begin{tikzpicture}[
                        baseline = (current bounding box.center),
                        thick,
                    ]
                        \draw (0,0) -- (1,0);
                        \node[left] at (0,0) {$3$};
                    \end{tikzpicture}
                }}
                \boxtimes 
                \vcenter{\hbox{
                    \begin{tikzpicture}[
    baseline = (current bounding box.center),
    scale = 0.5,
    thick,
]
    \draw (4,1) ..controls +(180:0.5) and +(90:0.5).. (3,0);
    \draw (3,0) ..controls +(-90:0.5) and +(180:0.5).. (4,-1);
    \draw (2,0) -- (3,0);

    \node[left] at (2,0) {$4$};
    \node[right] at (4,1) {$2$};
    \node[right] at (4,-1) {$2$};
\end{tikzpicture}
                }}
                \right)\\
                &\star\left(
                    \vcenter{\hbox{
                        \begin{tikzpicture}[
    baseline = (current bounding box.center),
    scale = 0.5,
    thick,
]
    \draw (0,1) ..controls +(0:0.5) and +(90:0.5).. (1,0);
    \draw (1,0) ..controls +(-90:0.5) and +(0:0.5).. (0,-1);
    \draw (2,0) -- (1,0);

    \node[left] at (0,1) {${3}$};
    \node[left] at (0,-1) {$2$};
    \node[right] at (2,0) {$5$};
\end{tikzpicture}
                    }}
                    \boxtimes
                    \vcenter{\hbox{
                        \begin{tikzpicture}[
                            baseline = (current bounding box.center),
                            thick,
                        ]
                            \draw (0,0) -- (1,0);
                            \node[left] at (0,0) {$2$};
                        \end{tikzpicture}
                    }}
                \right)
        \end{align*}
    \end{example}
\subsection{Rickard complex}
In this section, we recall the construction of the Rickard complex of singular Soergel bimodules associated to a colored braid. To begin with, let $\mathrm{Br}_n$ denote the $n$-strand braid group. We first extend braid group to colored braid groupoid.

\begin{definition}[Colored braid groupoid]
    Let $S$ be a color set. Define the \emph{$S$-colored groupoid} $\mathfrak{Br}(S)$ as the category with the following data: 
    \begin{itemize}
        \item objects in $\mathfrak{Br}(S)$ are sequences of the form $\mathbf{a}=(a_1,a_2,\ldots,a_m)$, in which  $m\geq 1$ and $a_i\in S$ for every $1\leq i\leq m$;
        \item morphisms from $\mathbf{a}$ to $\mathbf{b}$ in $\mathfrak{Br}(S)$ are $S$-colored braids with right boundary $\mathbf{a}$ and left boundary $\mathbf{b}$. To be precise, 
        \[
            \mathrm{Hom}_{\mathfrak{Br}(S)}(\mathbf{a},\mathbf{b}):=\{\beta\in\mathrm{Br}_m: a_i=b_{\beta(i)},\forall 1\leq i\leq m\}
        \]
        where $\mathbf{a}=(a_1,a_2,\ldots,a_m)$ and $\mathbf{b}=(b_1,b_2,\ldots,b_m)$. 
    \end{itemize}
\end{definition}
The colored braid groupoid $\mathfrak{Br}(S)$ is generated by the following \emph{Artin generators} $\beta_i$'s, defined as 
\[
    \beta_i: (a_1,\ldots,a_{i-1},a_i,a_{i+1},a_{i+2},\ldots,a_m) \mapsto (a_1,\ldots,a_{i-1},a_{i+1},a_{i},a_{i+2},\ldots,a_m).
\]
A \emph{colored braid word} of a colored braid $\beta$ is a decomposition of colored braid into Artin generators $\beta=\beta_{i^1}^{\varepsilon_1}\cdots\beta_{i_r}^{\varepsilon_r}$, written as ${(\beta)}_{\mathbf{a}}$. From now on, we set $S=\mathbb{Z}_{\geq 1}$, so that the objects in $\mathfrak{Br}(S)$ are identical to the objects in the 2-category $\mathrm{SSBim}$. Let $\mathcal{C}(\mathrm{SSBim})$ denote the category of chain complexes of singular Soergel bimodules. 
\begin{definition}[1-strand Rickard complex]
    Let $a\geq 1$ be an positive integer. Define the Rickard complex corresponding to the 1-strand braid with coloring $a$ from left to right as 
    \[
        \left\llbracket\vcenter{\hbox{
            \begin{tikzpicture}[
                baseline = (current bounding box.center),
                thick,
            ]
                \draw (0,0) -- (2,0);
                \node[right] at (2,0) {${a}$};
            \end{tikzpicture}
        }}
        \right\rrbracket := \left(
            0 \to \mathrm{Sym}(\mathbb{X}) \to 0
        \right),
    \]
    where $\mathbb{X}=\{x_1,\ldots,x_a\}$ is an alphabet of size $a$, and the non-zero term is in homological degree $0$.
\end{definition}
\begin{definition}[2-strand Rickard complex]
    Let $a,b\geq 1$ be two integers. Define the Rickard complex $C_{a,b}$ corresponding to a positive or negative crossing as the following bounded complex in $\mathcal{C}(\mathrm{SSBim})$ 
    \begin{align*}
        C_{a,b}&=\left\llbracket 
        \vcenter
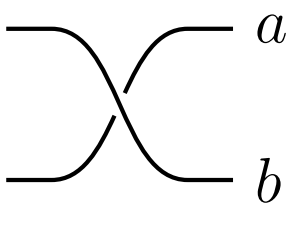
{\hbox{
            
\begin{tikzpicture}[
    baseline = (current bounding box.center),
    every braid/.style ={
        thick,
        rotate = 270,
        braid/anchor = center,
    }
]
    \pic (mybraid) at (0,0) {
        braid = {s_1^{-1}}
    };
    \node[right] at (mybraid-1-s) {${a}$};
    \node[right] at (mybraid-2-s) {${b}$};
\end{tikzpicture}
 
        }}
        \right\rrbracket
        :=\left(C_{a,b}^0
        \overset{\chi_0^+}{\to} \cdots 
        \overset{\chi_0^+}{\to} q^{-k}t^{k} C_{a,b}^k 
        \overset{\chi_0^+}{\to} q^{-k-1}t^{k+1} C_{a,b}^{k+1}
        \overset{\chi_0^+}{\to} \cdots
        \right)
    \end{align*}
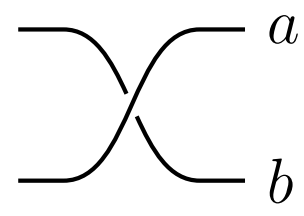
    \begin{align*}
        C_{a,b}^\vee &=
        \left\llbracket 
        \vcenter{\hbox{
            
\begin{tikzpicture}[
    baseline = (current bounding box.center),
    every braid/.style ={
        thick,
        rotate = 270,
        braid/anchor = center,
    }
]
    \pic (mybraid) at (0,0) {
        braid = {s_1}
    };
    \node[right] at (mybraid-1-s) {${a}$};
    \node[right] at (mybraid-2-s) {${b}$};
\end{tikzpicture}

        }}
        \right\rrbracket=\left(\cdots 
        \overset{\chi_0^-}{\to} q^{k+1}t^{-k-1} C_{a,b}^{k+1}
        \overset{\chi_0^-}{\to} q^{k}t^{-k} C_{a,b}^{k}
        \overset{\chi_0^-}{\to} \cdots
        \overset{\chi_0^-}{\to} C_{a,b}^0\right)
    \end{align*}
    where $C_{a,b}^k$ is the singular Soergel bimodule associated to the following web
    \[
        C_{a,b}^k:=
        \vcenter
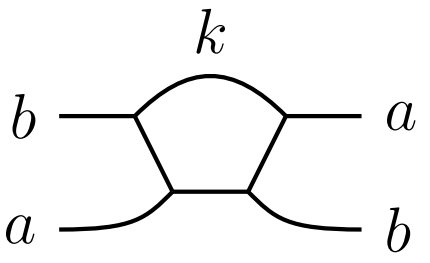
{\hbox{
                \begin{tikzpicture}[
        thick, 
        scale = 0.5, 
        baseline= (current bounding box.center)]
        \draw (0,0) -- (1,0);
        \draw (1,0) ..controls +(45:1) and +(135:1).. (3,0)
        node[midway, above] {$k$};
        \draw (3,0) -- (4,0);
        \draw (1,0) -- (1.5,-1);
        \draw (3,0) -- (2.5,-1);
        \draw (1.5,-1) -- (2.5,-1);
        \draw (0,-1.5) ..controls +(1,0) and +(225:0.5).. (1.5,-1);
        \draw (4,-1.5) ..controls +(-1,0) and +(315:0.5).. (2.5,-1);

        \node[left] at (0,0) {$b$};
        \node[left] at (0,-1.5) {$a$};
        \node[right] at (4,0) {$a$};
        \node[right] at (4,-1.5) {$b$};
    \end{tikzpicture}}
        }
    \]
    and the differentials in $C_{a,b}$ and $C_{a,b}^\vee$ are given by $\chi_0^\pm:C_{a,b}^k \to C_{a,b}^{k\pm 1}$, respectively, where $\chi_0^{\pm}$ are bimodule homomorphisms between singular Soergel bimodules $C_{a,b}^k$ defined in~\cite[Section 3.5]{hogancamp2024link}. Their precise form will not be used in the sequel. 
\end{definition}

\begin{example}
    The Rickard complex of a crossing with $(2,2)$-coloring is given by 
    \begin{align*}
        C_{2,2} 
        &=\left\llbracket
        \vcenter
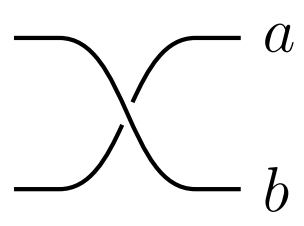
{\hbox{
            
\begin{tikzpicture}[
    baseline = (current bounding box.center),
    every braid/.style ={
        thick,
        rotate = 270,
        braid/anchor = center,
    }
]
    \pic (mybraid) at (0,0) {
        braid = {s_1^{-1}}
    };
    \node[right] at (mybraid-1-s) {${a}$};
    \node[right] at (mybraid-2-s) {${b}$};
\end{tikzpicture}
 
        }}\right\rrbracket\\
        &=\left(
        \vcenter
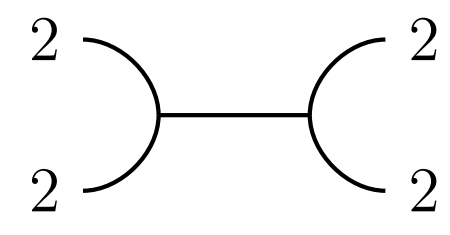
{\hbox{
            
\begin{tikzpicture}[
    baseline = (current bounding box.center),
    thick,
    scale = 0.5,
]
    \draw (0,1) ..controls +(0:0.5) and +(90:0.5).. (1,0);
    \draw (1,0) ..controls +(-90:0.5) and +(0:0.5).. (0,-1);
    \draw (4,1) ..controls +(180:0.5) and +(90:0.5).. (3,0);
    \draw (3,0) ..controls +(-90:0.5) and +(180:0.5).. (4,-1);
    \draw (1,0) -- (3,0);

    \node[left] at (0,1) {${2}$};
    \node[left] at (0,-1) {${2}$};
    \node[right] at (4,1) {${2}$};
    \node[right] at (4,-1) {${2}$};
\end{tikzpicture}

        }}
        \overset{\chi_0^+}{\to}
        q^{-1}t
        \vcenter
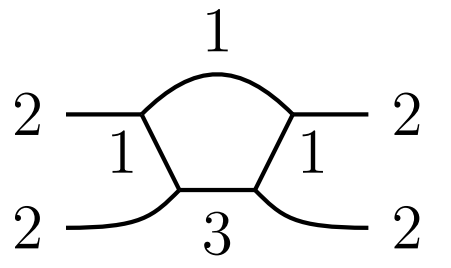
{\hbox{
            
    \begin{tikzpicture}[thick, scale = 0.5]
        \draw (0,0) -- (1,0);
        \draw (1,0) ..controls +(45:1) and +(135:1).. (3,0)
        node[midway, above] {$1$};
        \draw (3,0) -- (4,0);
        \draw (1,0) -- (1.5,-1)
        node[midway, left] {$1$};
        \draw (3,0) -- (2.5,-1)
        node[midway, right] {$1$};
        \draw (1.5,-1) -- (2.5,-1)
        node[midway, below] {$3$};
        \draw (0,-1.5) ..controls +(1,0) and +(225:0.5).. (1.5,-1);
        \draw (4,-1.5) ..controls +(-1,0) and +(315:0.5).. (2.5,-1);

        \node[left] at (0,0) {$2$};
        \node[left] at (0,-1.5) {$2$};
        \node[right] at (4,0) {$2$};
        \node[right] at (4,-1.5) {$2$};
    \end{tikzpicture}

        }}
        \overset{\chi_0^+}{\to}
        q^{-2}t^2
        \vcenter
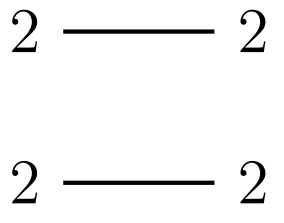
{\hbox{
            
    \begin{tikzpicture}[
        thick, 
        baseline = (current bounding box.center),
        scale = 0.5,
        ]
        \draw (0,1) -- (2,1)
        node[pos = 0, left] {$2$}
        node[pos = 1, right] {$2$};
        \draw (0,-1) -- (2,-1)
        node[pos = 0, left] {$2$}
        node[pos = 1, right] {$2$};
    \end{tikzpicture}

        }}
        \right),
    \end{align*}
    where $q$ is the internal grading shift, and $t$ is the homological grading shift.
\end{example}
\begin{remark}
    The complex $C_{a,b}$ is clearly a bounded complex, as the range of $k$ is $0\leq k\leq \min(a,b)$. In addition, the Soergel bimodule $C_{a,b}^k$ is also called a \emph{ladder web}, which can be described as the quotient ring of symmetric polynomial algebras, as for example in~\cite[Section 3.3]{conners2025doubly}. 
\end{remark}
Note that every braid can be represented as a composition of Artin generators, which are vertical stacking of the 2-strand crossings and the 1-strand identity braids defined above. Thus the Rickard complex of a braid can be defined as the tensor product of Rickard complexes of the Artin generators.
\begin{definition}[Rickard complex of a braid]\label{def:rickard-complex}
    Let $\beta$ be an $m$-strand colored braid, and write 
    \[
        \beta = \beta_{i_1}^{\varepsilon_1} \cdot \beta_{i_2}^{\varepsilon_2}\cdots \beta_{i_r}^{\varepsilon_r},
    \]
    where $\beta_{i_j}$ is an $m$-strand braid with a single positive crossing exchanging the $i_j$-th and $(i_j+1)$-th strands for every $1\leq j\leq r$, and $\varepsilon_1,\ldots,\varepsilon_r\in\{\pm 1\}$. Define the \emph{Rickard complex} of $\beta$ as
    \[
        C(\beta_{\mathbf{a}}):=\llbracket \beta_{i_1}^{\varepsilon_1} \rrbracket \star \llbracket \beta_{i_2}^{\varepsilon_2}\rrbracket \star\cdots \star \llbracket \beta_{i_r}^{\varepsilon_r} \rrbracket,
    \]
    where $\star$ is the horizontal composition of singular Soergel bimodules defined in Definition~\ref{def:singular-soergel-bimodule}, which is the tensor product over the action of the common symmetric polynomial algebra.
\end{definition}
\subsection{Hochschild homology}\label{section:hochschild-homology}
    In this section, we define a functor $\mathrm{HH}_*$ called \emph{Hochschild homology} from the category of graded $R$-bimodules to the category of triply graded vector spaces over $\mathbb{C}$, where $R$ is a commutative graded algebra over $\mathbb{C}$. In particular, Hochschild homology for singular Soergel bimodules (or say modules over $R\otimes R^{\mathrm{op}}$) can be computed using Koszul factorization. Moreover, there is a reduced version of Hochschild homology functor $\mathrm{HH}_{\neg p}$ for reduced version of colored HOMFLY homology.


    Note that every $R$-bimodule can also be regarded as a $R\otimes R^{\mathrm{op}}$-module, and thus one can define a graded $R$-bimodule $M$ its Hochschild homology $\mathrm{HH}_*(M)$ as the following bigraded vector spaces.
    \begin{definition}[Hochschild homology]
        Let $R$ be a graded commutative algebra over $\mathbb{C}$, and $M$ be a graded $R$-bimodule. The \emph{$i$-th Hochschild homology} of $M$ is defined as a $\mathbb{Z}_q$-graded vector space 
        \[
            \mathrm{HH}_{i}(R, M):=\mathrm{Tor}_i^{R\otimes R^{\mathrm{op}}}(R,M),
        \]
        where the graded $R$-bimodule structure of $R$ is given by the left and right multiplication. Moreover, the \emph{Hochschild homology} of $M$ is defined as a $\mathbb{Z}_a\times\mathbb{Z}_q$-bigraded vector space
        \[
            \mathrm{HH}_*(R, M):=\bigoplus_{i\in\mathbb{Z}} a^i\mathrm{HH}_i(R, M).
        \]
    \end{definition}

    \begin{remark}
        We will use $\mathrm{HH}_*(M)$ to denote $\mathrm{HH}_*(R, M)$, if $M$ is a graded $R$-bimodule and the graded algebra $R$ is clear from the context.
    \end{remark}
    Hochschild homology can be expressed in terms of the homology of a \emph{Koszul complex}. This plays an important role in the extension of Hochschild homology to other versions, especially the reduced version $\mathrm{HH}_{\neg p}$ at a point $p$ on the braid. 
    \begin{definition}[Koszul complex]
        Let $\mathbf{x}=\{x_1,\ldots,x_n\}$ be a set of $\mathbb{Z}_q$-homogeneous elements in $R$, with $\deg_q(x_i)=d_i$ for every $1\leq i\leq n$.  
        The \emph{Hochschild Koszul complex} $\mathcal{Z}^a_\mathbf{x}$ is defined as the $\mathbb{Z}_a$-graded complex 
        \[
            \mathcal{Z}^a_{\mathbf{x}}:=\sideset{}{_R}\bigotimes_{i=1}^n\left( R \overset{\cdot x_i}{\rightarrow} aq^{-d_i}R\right)
        \]
        of $\mathbb{Z}_q$-graded $R$-algebras, where $a$ is the Hochschild grading. 
    \end{definition}
    
    For a colored braid $\beta$ with initial object $(a_1,\ldots,a_m)$, let $R_\beta:=\mathrm{Sym}(\mathbb{X}_1\mid\cdots\mid\mathbb{X}_m)$ be the symmetric polynomial algebra associated to $\beta$, where $\mathbb{X}_i$ is an alphabet of size $|\mathbb{X}_i|=a_i$ for every $1\leq i\leq m$. Moreover, we identify $R_\beta\otimes R_\beta \cong \allowbreak\mathrm{Sym}(\mathbb{X}_1\mid\cdots\mid\mathbb{X}_m\mid\mathbb{X}'_1\mid\cdots\mid\mathbb{X}'_m)$ with the symmetric polynomial algebra where $\mathbb{X}'_i$ is also an alphabet $|\mathbb{X}'_i|=a_i$ for every $1\leq i\leq m$. Then the Hochschild homology of an $(R_\beta,R_\beta)$-bimodule $M$ can be expressed in terms of the homology of Koszul complex $\mathcal{Z}_\mathbf{x}$ for certain $\mathbf{x}$, as stated in the following lemma.

    \begin{lemma}[{\cite[Lemma 3.7]{wedrich2019exponential}}]\label{lem:hh-koszul}
        Let $\mathbf{x}=\{e_j(\mathbb{X}_i)-e_j(\mathbb{X}'_i)\}\subset R_\beta\otimes R_\beta$ be a set of symmetric polynomials of quantum degree $2j$, where $1\leq j\leq a_i$, $1\leq i\leq m$, and $e_j(\mathbb{X})$ is the $j$-th elementary symmetric polynomial in $\mathbb{X}$. Let $M$ be a graded $R_\beta$-bimodule. Then there is an isomorphism between bigraded vector spaces
        \[
            \mathrm{HH}_*(R_\beta, M)\cong \mathrm{H}_*(M\otimes_{R_\beta} \mathcal{Z}^a_\mathbf{x}).
        \]
    \end{lemma}

    The reduced version of Hochschild homology is defined on a $m$-strand colored braid $\beta$ with a marked point $p$ on it. For simplicity, we use $\mathbb{X}_p$ to denote the alphabet associated to the strand where the marked point $p$ lies, and similarly denote $a_p:=|\mathbb{X}_p|$. 

    \begin{example}
        Let $\beta$ be a colored braid with a marked point $p$ lying on the first strand, as drawn in red in Figure~\ref{fig:example-for-braid-with-marked-point}. In this case, we write $\mathbb{X}_p=\mathbb{X}_1$ and $a_p=1$. 
        \begin{figure}[h]
            \centering 
            \includegraphics[width = 0.3\textwidth]{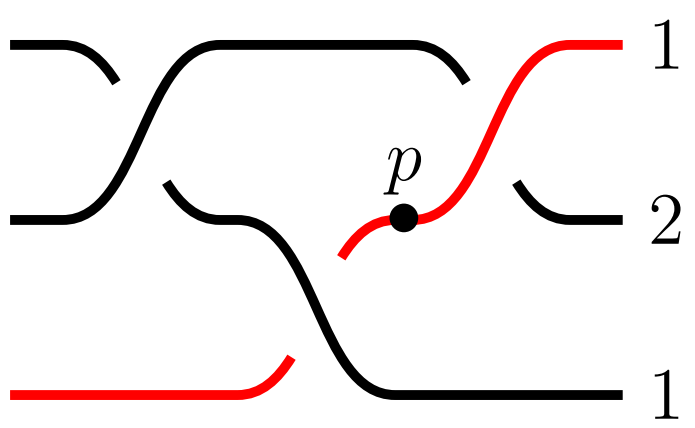}
            \caption{Braid with a marked point}\label{fig:example-for-braid-with-marked-point}
        \end{figure}
    \end{example}
    \begin{definition}[Reduced Hochschild homology]
        Let $\beta$ be a colored braid with a marked point $p$ away from all crossings and with coloring $(a_1,\ldots,a_m)$. Set $\hat{\mathbf{x}}:=\{e_j(\mathbb{X}_i)-e_j(\mathbb{X}'_i): \mathbb{X}_i\neq\mathbb{X}_p, 1\leq j\leq a_i\}$. The \emph{reduced Hochschild homology} of an $R_\beta$-bimodule $M$ at $p$ is defined as the $\mathbb{Z}_a\times\mathbb{Z}_q$-graded vector space
        \[
            \mathrm{HH}_{\neg p}(M):=\mathrm{H}_*(M\otimes_{R_\beta} \mathcal{Z}^a_{\hat{\mathbf{x}}})
        \]
        obtained by taking the homology of the tensor product of $M$ and the Koszul complex $\mathcal{Z}^a_{\hat{\mathbf{x}}}$.
    \end{definition}
    The reduced Hochschild homology functor $\mathrm{HH}_{\neg p}$ represents a graphical operation closing the braid at all ends except at the end associated to $p$. To be precise, if $p$ lies on the $i$-th strand, then we connect the $j$-th end on the left and the $j$-th end on the right, for every $j\neq i$. See the following example.  
    \begin{example}[Graphical calculus for reduced Hochschild functor]
        Suppose that $\beta$ is a colored braid with three strands as in Figure~\ref{fig:example-for-braid-with-marked-point}. Then $\mathrm{HH}_{\neg p}(\beta)$ is represented by Figure~\ref{fig:closure-except}. 
        \begin{figure}[h]
            \centering 
            \includegraphics[width = 0.5\textwidth]{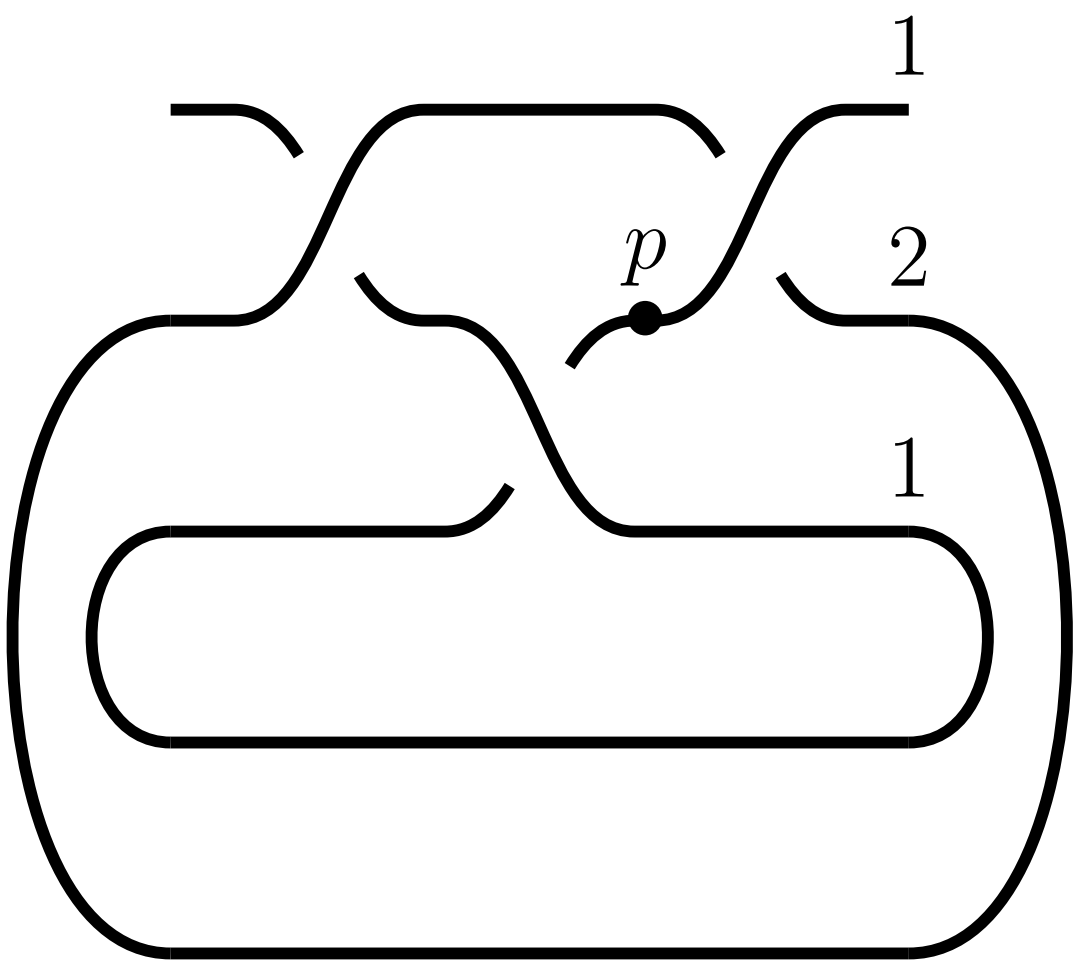}
            \caption{$\mathrm{HH}_{\neg p}(\beta)$}\label{fig:closure-except}
        \end{figure}
        \FloatBarrier{}
    \end{example}

\subsection{Colored HOMFLY link homology}\label{section:colored-homfly-homology}
In this section, we recall the construction of the triply graded colored link homology introduced in~\cite{webster2017geometric} following~\cite{conners2025doubly}. Let us first introduce some notations following~\cite{conners2025doubly}. 

\begin{definition}[Balanced braid]
    Let ${}_{\mathbf{c}}\beta_{\mathbf{b}}$ denote a colored braid with left boundary color $\mathbf{c}$ and right boundary color $\mathbf{b}$. Then ${}_{\mathbf{c}}\beta_\mathbf{b}$ is balanced if $\mathbf{c}=\mathbf{b}$.
\end{definition}
\begin{example}
    Figure~\ref{fig:4-strand-balanced-braid} is an example for a balanced braid, with $\mathbf{b}=(1,2,1,2)$ and $\mathbf{c}= (1,2,1,2)$, whose closure is a link with two link components. 
    \begin{figure}[H]
        \centering 
        \includegraphics[width = 0.3\textwidth]{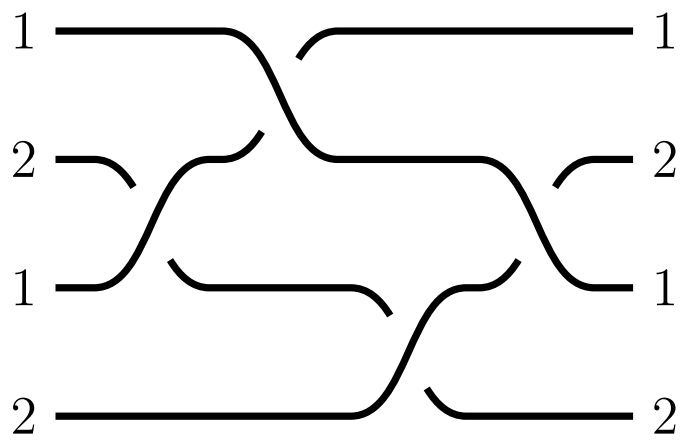}
        \caption{4-strand balanced colored braid}\label{fig:4-strand-balanced-braid}
    \end{figure}
    \FloatBarrier{}
\end{example}

Note that only balanced braids have well-defined closure. For simplicity, a balanced colored braid ${}_{\mathbf{c}}\beta_{\mathbf{b}}$ is also written as $\beta_{\mathbf{b}}$. Suppose that the closure of a balanced braid $\beta$ is a colored link $\mathcal{L}$. Then strands in $\beta$ corresponding to the same link component in $\mathcal{L}$ must have the same coloring. To better describe behavior of links in colored braid, we introduce the following notion.

\begin{definition}
    Let $\omega\in\mathfrak{S}_m$ be a permutation. Define an equivalence relation $\sim$ on the set $\{1,2,\ldots,m\}$ as $i\sim j$ if and only if there exists some $k\geq 0$, such that $i=\omega^k(j)$. Moreover, write 
    \[
        \Omega(\omega):=\{1,2,\ldots,m\}/\sim
    \]
    as the set of cycles. If $\mathbf{b}=(b_1,b_2,\ldots,b_m)$ satisfies $\omega(\mathbf{b})=\mathbf{b}$, namely $b_{\omega(i)}=b_i$ holds for every $1\leq i\leq m$, then write $b_{[i]}:=b_i$.
\end{definition}
\begin{remark}
    In particular, if $\omega$ is the permutation corresponding to the braid $\beta$ whose closure is a link $\mathcal{L}$, then there is a one-to-one correspondence between $\Omega(\omega)$ and the set of link components $\pi_0(\mathcal{L})$. If $\mathcal{L}$ is in addition colored, then the boundary object $\mathbf{b}$ satisfies $\omega(\mathbf{b})=\mathbf{b}$. 
\end{remark}

We introduce some numerical invariants of colored braids following~\cite{hogancamp2024link}.

\begin{definition}
    Let $\beta_{\mathbf{b}}$ be an $m$-strand colored braid, with boundary object $\mathbf{b}=(b_1,b_2,\ldots,b_m)$ on the right. Define the colored writhe $\varepsilon(\beta)$ of $\beta_\mathbf{b}$ as 
    \[
        \varepsilon(\beta_\mathbf{b}):=\sum_{c\in\mathrm{cr}(\beta)} s(c)a(c)b(c),
    \]
    where $\mathrm{cr}(\beta)$ is the set of all crossings in the braid $\beta$, $s(c)\in\{\pm 1\}$ is the sign of the crossing $c$, and $\{a(c),b(c)\}$ is the set of the colors of the two strands that intersect at $c$. If $\beta$ is in addition balanced, then we can define 
    \[
        n(\beta)=\sum_{[i]\in\Omega(\omega)} b_{[i]},
    \]
    \[
        N(\beta):=\sum_{1\leq i\leq m} b_i,
    \]
    and 
    \[
        Q(\beta):=\sum_{1\leq i\leq m}^{} b_i(b_i+1).
    \]
\end{definition}
With the above numerical invariants and Rickard complexes of colored braids, we are now in a position to define colored HOMFLY homology of a link.
\begin{definition}[Colored HOMFLY homology]
    Let $\mathcal{L}$ be a colored link that can be written $\mathcal{L}=\hat{\beta_{\mathbf{b}}}$ as the closure of a balanced braid $\beta$. Define the colored Khovanov-Rozansky complex of $\beta$ as 
    \begin{equation}
        C_{\mathrm{KR}}(\beta):={(at^{-1})}^{\frac{1}{2}\left(\varepsilon(\beta)+N(\beta)-n(\beta)\right)} q^{-\varepsilon(\beta)}\cdot\mathrm{HH}(C(\beta)),
    \end{equation}
    where $C(\beta)$ is the Rickard complex of the braid $\beta$, and the colored Khovanov-Rozansky homology as the homology of the above complex, namely
    \begin{equation}
        H_{\mathrm{KR}}(\beta):=\mathrm{H}(C_{\mathrm{KR}}(\beta)).
    \end{equation}
\end{definition}
The invariance of the choice of braid presentation for links is proved by Webster and Williamson in~\cite{webster2017geometric} as stated in the theorem below. 
\begin{theorem}[{\cite[Theorem 1.1]{webster2017geometric}}]
    The colored HOMFLY homology is a link invariant, up to isomorphisms of triply graded vector spaces, and its graded Euler characteristic is the colored HOMFLY polynomial.
\end{theorem}
\begin{remark}
When every link component of a link $\mathcal{L}$ is colored by $1$, the column-colored homology defined above recovers the original uncolored homology introduced by Khovanov and Rozansky in~\cite{khovanov2008matrix}. 
\end{remark}
Having introduced the notion of colored HOMFLY homology, we now focus on the reduced version of colored HOMFLY homology.
\begin{definition}[Homological Koszul complex]
    Let $\beta$ be a colored braid with a marked point $p$ on it, and $R:=R_\beta$ be the symmetric polynomial algebra associated to $\beta$. Let $\mathbf{x}_p:={\{e_i(\mathbb{X}_p)\}}_{1\leq i\leq a_p}$ be the set of all elementary symmetric polynomials in $\mathbb{X}_p$. Define the \emph{homological Koszul complex} $\mathcal{Z}^t_{\mathbf{x}_p}$ as the $\mathbb{Z}_t$-graded complex
    \[
        \mathcal{Z}^t_{\mathbf{x}_p}:=\sideset{}{_{R}} \bigotimes_{i=1}^{a_p} \left(
            R \overset{\cdot e_i(\mathbb{X}_p)}{\longrightarrow} tq^{-2i} R
        \right),
    \]
    of $\mathbb{Z}_q$-graded $R$-bimodules.
\end{definition}
\begin{definition}[Reduced column-colored Rickard complex]
    Let $\beta$ be a colored braid with a marked point $p$, and $\mathbf{x}_p:={\{e_i(\mathbb{X}_p)\}}_{1\leq i\leq a_p}$ be the set of all elementary symmetric polynomials in $\mathbb{X}_p$. Define the \emph{reduced column-colored Rickard complex} $C_p(\beta)$ of $\beta$ at $p$ as 
    \[
        C_p(\beta):= C(\beta)\otimes_R \mathcal{Z}^t_{\mathbf{x}_p},
    \]
    where $C(\beta)$ is the Rickard complex associated to $\beta$. Moreover, the \emph{reduced colored Khovanov-Rozansky complex} of $\beta$ at $p$ is defined as the $\mathbb{Z}_t$-graded complex of $\mathbb{Z}_a\times\mathbb{Z}_q$-graded vector spaces
    \[
        C_{\mathrm{KR}}(\beta(p)):={(at^{-1})}^{\frac{1}{2}\left(\varepsilon(\beta)+N(\beta)-n(\beta)\right)} q^{-\varepsilon(\beta)}\cdot\mathrm{HH}_{\neg p}(C_p(\beta)),
    \]
    and the \emph{reduced colored Khovanov-Rozansky homology} of $\beta$ at $p$ is defined as the homology of the above complex with respect to the homological grading, namely
    \[
        H_{\mathrm{KR}}(\beta(p)):=\mathrm{H}(C_{\mathrm{KR}}(\beta(p))).
    \]
\end{definition}
This reduced version of colored Khovanov-Rozansky homology is a link invariant as stated in the following theorem~{\cite[Theorem 3.44]{conners2025doubly}}.
\begin{theorem}[Reduced column-colored HOMFLY homology]
    Let $\mathcal{L}$ be a colored link with a marked point $p$ on it, and $\beta$ be a colored braid with closure $\hat{\beta}=\mathcal{L}$. Then $H_{\mathrm{KR}}(\beta(p))$ is an invariant of colored link $\mathcal{L}$ and the link component $\mathcal{L}_p\subset\mathcal{L}$ of $p$, up to isomorphisms of $\mathbb{Z}_a\times\mathbb{Z}_q\times\mathbb{Z}_t$-graded vector spaces. Therefore, we can define the \emph{reduced column-colored HOMFLY homology} of $\mathcal{L}$ at $p$ as
    \[
        \overline{\mathrm{HHH}}(\mathcal{L}(p)):= H_{\mathrm{KR}}(\beta(p)),
    \]
    where $\beta$ is any colored braid with closure $\hat{\beta}=\mathcal{L}$.
\end{theorem}
\begin{remark}
    For a link $L$ with $r$ components with a marked point $p$ on it, we denote by
    \[
        \overline{\mathrm{HHH}}_{\bigwedge^{a_1},\bigwedge^{a_2},\ldots,\bigwedge^{a_r}}(L(p))
    \] 
    the reduced column-colored HOMFLY homology of $L$ at $p$, where the $i$-th link component carries color $a_i$ for every $1\leq i\leq r$. In particular, if $L=K$ is a knot, then the marked point $p$ is always omitted and write 
    \[
        \overline{\mathrm{HHH}}_{\bigwedge^{a}}(K):=\overline{\mathrm{HHH}}_{\bigwedge^{a}}(K(p)),
    \]
    as there is only one choice for link component. 
\end{remark}

\subsection{Quantum multinomial formula}\label{section:quantum-multinomial-formula}
We closed this section by establishing a quantum multinomial formula, which will be helpful in our main computation. This quantum multinomial formula is likely wellknown, but we reproduce a proof here. To begin with, we define the notion of strict inversion and inversion number of a sequence. 
\begin{definition}[Inversion]
    Let $\Delta:\{1,2,\ldots,m\}\to \{0,1,\ldots,p\}$ be a sequence of length $m$. A pair $(i,j)$ is called a \emph{strict inversion} of $\Delta$ if $i<j$ and $\Delta(i)>\Delta(j)$. The \emph{strict inversion number} of $\Delta$ is defined as the number of strict inversions in $\Delta$, denoted by 
    \[
        \mathrm{inv}(\Delta):=\#\{(i,j): i<j, \Delta(i)>\Delta(j)\}.
    \]
    Similarly, a pair $(i,j)$ is called a \emph{inversion} of $\Delta$ if $i<j$ and $\Delta(i)\geq \Delta(j)$. The \emph{inversion number} of $\Delta$ is defined as the number of weak inversions in $\Delta$, denoted by
    \[
        \mathrm{Inv}(\Delta):=\#\{(i,j): i<j, \Delta(i)\geq \Delta(j)\}.
    \]
\end{definition}
\begin{lemma}[Quantum multinomial theorem]\label{lem:multinomial-theorem}
    Let $k_0=m\geq k_1\geq k_2\geq\cdots\geq k_p\geq 0=k_{p+1}$, and $m\geq 0$. Suppose that $\Delta :\{1,2,\ldots,m\}\to \{0,1,\ldots,p\}$ is a sequence of length $m$, such that 
    \[
        \#\Delta^{-1}(i) = k_i - k_{i+1}, \quad\forall 0\leq i\leq p.
    \]
    For fixed $\kappa=(k_0,k_1,k_2,\ldots,k_p)$, write $S_\kappa$ as the set of all $\Delta$ of length $k_0$ with $\# \Delta^{-1}(i)=k_i -k_{i+1}$, for every $0\leq i\leq p$. Then 
    \begin{equation}
        \sum_{\Delta\in S_\kappa} Q^{\mathrm{inv}(\Delta)} = \sqbinom{k_0}{k_1}\sqbinom{k_1}{k_2}\cdots\sqbinom{k_{p-1}}{k_p}.\label{eq:multinomial}
    \end{equation}
\end{lemma}
\begin{proof}
    We prove the lemma by induction on $p$. For $p=0$, there is exactly one $\Delta$, with $\mathrm{inv}(\Delta)=0$. Then 
    \begin{equation}
        \mathrm{LHS}=Q^{0}=1=\sqbinom{0}{0}=\mathrm{RHS}.\nonumber
    \end{equation}
    For $p=1$, the sequence $\Delta:\{1,2,\ldots,m\}\to \{0,1\}$, with $k_1\leq m=k_0$ has exactly $k_1$ ones, and $(m-k_1)$ zeros in the sequence. Define 
    \[
        F(m,k_1):=\sum_{\Delta\in S_{m,k_1}}Q^{\mathrm{inv(\Delta)}}.
    \]
    Clearly, (\ref{eq:multinomial}) holds for $(m,k_1)=(1,0),(1,1)$. Note that 
    \begin{align}
        F(m,k_1) 
        &= \sum_{\Delta(m)=0} Q^{\mathrm{inv}(\Delta')}\cdot Q^{k_1} + \sum_{\Delta(m)=1} Q^{\mathrm{inv}(\Delta')}\nonumber\\
        &=Q^{k_1} F(m-1,k_1) + F(m-1,k_1-1), \label{eq:F(m,k)}
    \end{align}   
    where $\Delta'$ is obtained from $\Delta$ by removing the last digit. Moreover, 
    \begin{align}
        \sqbinom{m}{k_1} 
        &=\frac{[m]!}{[k_1]![m-k_1]!} \nonumber\\
        &=\frac{[m]}{[k_1]}\cdot\frac{[m-1]!}{[k_1-1]![m-k_1]!}\nonumber\\
        &=\frac{1-Q^m}{1-Q^{k_1}}\cdot\frac{[m-1]!}{[k_1-1]![m-k_1]!}\nonumber\\
        &=\frac{1-Q^{k_1}+Q^{k_1}-Q^{m}}{1-Q^{k_1}}\cdot \frac{[m-1]!}{[k_1-1]![m-k_1]!}\nonumber\\
        &=\sqbinom{m-1}{k_1-1} + \frac{Q^{k_1}(1-Q^{m-k_1})}{1-Q^{k_1}}\sqbinom{m-1}{k_1-1}\nonumber\\
        &=\sqbinom{m-1}{k_1-1} + Q^{k_1}\frac{[m-1]![m-k_1]}{[k_1-1]![m-k_1]![k_1]}\nonumber\\
        &=\sqbinom{m-1}{k_1-1} + Q^{k_1} \sqbinom{m-1}{k_1}. \label{eq:binomial-m-k}
    \end{align}
    By comparing (\ref{eq:F(m,k)}) and (\ref{eq:binomial-m-k}), we have $F(m,k_1)=\sqbinom{m}{k_1}$, which indicates that (\ref{eq:multinomial}) holds for $p=1$. Assume that (\ref{eq:multinomial}) holds for $0,1,\ldots,p-1$. Then for $\Delta:\{1,2,\ldots,m\}\to\{0,1,\ldots,p\}$, we have 
    \begin{align}
        &~~\mathrm{inv}(\Delta) \nonumber\\
        &=\#\{(i,j): i<j,\Delta(i)>\Delta(j)\}\nonumber\\
        &=\#\{(i,j): i<j,\Delta(i)>\Delta(j),\Delta(i)\neq p\} + \#\{(i,j): i<j,\Delta(i)>\Delta(j), \Delta(i)=p\}\nonumber\\
        &=\mathrm{inv}(\Delta\backslash\Delta^{-1}(p)) + \#\{(i,j):i<j, i\in\Delta^{-1}(p), j\notin\Delta^{-1}(p)\},\label{eq:inversion-of-delta-p-m}
    \end{align}
    where $\Delta\backslash\Delta^{-1}(p)$ is the sequence obtained from $\Delta$ by removing all the position $i$ which takes value $p$ in $\Delta$. Define $\tilde{\Delta} : \{1,2,\ldots,m\} \to \{0,1\}$ by
    \begin{equation*}
        \tilde{\Delta}(i):=
        \begin{cases}
            0, & \Delta(i)\leq p-1\\
            1, & \Delta(i)= p
        \end{cases}.
    \end{equation*}
    Then, 
    \begin{align}
        \sum_{\Delta\in S_\kappa} Q^{\mathrm{inv}(\Delta)} 
        &\overset{(\ref{eq:inversion-of-delta-p-m})}{=}
        \sum_{\Delta\in S_\kappa} Q^{\mathrm{inv}(\Delta\backslash\Delta^{-1}(p))}\cdot Q^{\#\{(i,j):i<j, i\in\Delta^{-1}(p), j\notin\Delta^{-1}(p)\}}\nonumber\\
        &=\sum_{\Delta\in S_\kappa} Q^{\mathrm{inv}(\Delta\backslash\Delta^{-1}(p))}\cdot Q^{\#\{(i,j):i<j, i\in\Delta^{-1}(p), j\notin\Delta^{-1}(p)\}}\nonumber\\
        &=\sum_{\Delta'\in S_{\kappa'}}\sum_{\tilde{\Delta}\in S_{m,k_p}} Q^{\mathrm{inv}(\Delta')} \cdot Q^{\mathrm{inv}(\tilde{\Delta})}\nonumber\\
        &=\sum_{\Delta'\in S_{\kappa'}} Q^{\mathrm{inv}(\Delta')}\sum_{\tilde{\Delta}\in S_{m,k_p}} Q^{\mathrm{inv}(\tilde{\Delta})},\nonumber
    \end{align}
    where $\kappa'=m-k_p=k'_0\geq k'_1\geq\cdots\geq k'_{p-1}\geq 0=k_p$ satisfies 
    \[
        \left\{\begin{aligned}
            k'_i-k'_{i+1}&=\#{\Delta'}^{-1}(i) = k_i-k_{i+1},\quad\forall 0\leq i\leq p-1\\
            k'_0&=m-k_p    \\
        \end{aligned}\right.
    \]
    Therefore, $k'_i=k_i-k_p$, for every $0\leq i\leq p-1$. By induction, we have 
    \begin{align*}
        \sum_{\Delta\in S_\kappa} Q^{\mathrm{inv}(\Delta)} 
        &=\left(\sqbinom{k'_0}{k'_1}\sqbinom{k'_1}{k'_2}\cdots\sqbinom{k'_{p-2}}{k'_{p-1}}\right)\cdot \sqbinom{m}{k_p}\\
        &=\left(\sqbinom{m-k_p}{k_1-k_p}\sqbinom{k_1-k_p}{k_2-k_p}\cdots\sqbinom{k_{p-2}-k_p}{k_{p-1}-k_p}\right)\cdot \sqbinom{m}{k_p}\\
        &=\left(\prod_{i=0}^{p-2}\frac{[k_i-k_p]!}{[k_i-k_{i+1}]![k_{i+1}-k_p]!}\right)\cdot \frac{[m]!}{[k_p]![m-k_p]!}\\
        &=\frac{[m]!}{[m-k_1]![k_1-k_2]!\cdots [k_{p-1}-k_p]! [k_p]!}\\
        &=\sqbinom{m}{k_1}\sqbinom{k_1}{k_2}\cdots\sqbinom{k_{p-1}}{k_p}.
    \end{align*}
    This completes the proof of equation (\ref{eq:multinomial}).
\end{proof}

\section{Computation of reduced colored homology}\label{chapter:computation}

\subsection{Recursive formula}
Before stating the recursive formulae, we introduce some notions on braids and permutations that will be useful throughout the recursive formulae and the computation.
\begin{definition}
    Let $v\in{\{0,1\}}^{l+m}$ be a $(01)$-sequence of length $l+m$, with $l$ zeros and $m$ ones. Suppose that the $i$-th zero in $v$ is the $a_i$-th digit, and the $j$-th one in $v$ is the $b_j$-th digit, for every $1\leq i\leq l$ and $1\leq j\leq m$. Then we define a permutation $\pi_v\in\mathfrak{S}_{l+m}$ by
    \[
        \pi_v(i)=a_i, \quad \pi_v(l+j)=b_j,
    \]
    for every $1\leq i\leq l$ and $1\leq j\leq m$. 
\end{definition}
\begin{example}
    Let $v=(10010)$, where $(a_1,a_2,a_3)=(2,3,5)$ and $(b_1,b_2)=(1,4)$. Then
   \[
        \pi_v =
        \vcenter
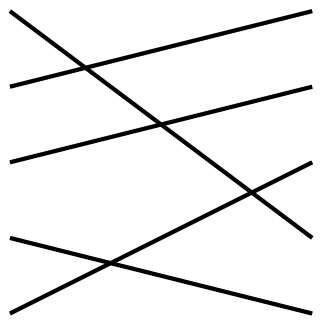
{\hbox{
        \begin{tikzpicture}[
            thick,
        ]
            \draw (0,0)--(2,-1.5);
            \draw (0,-0.5)--(2,0);
            \draw (0,-1)--(2,-0.5);
            \draw (0,-1.5)--(2,-2);
            \draw (0,-2)--(2,-1);
        \end{tikzpicture}
        }}
        =(12354)\in\mathfrak{S}_5,
    \]
    where the permutation $\pi_v$ is presented as a diagram from right to left.
\end{example}
\begin{definition}[Trace]
    Let $n\geq 1$, and $\omega\in\mathfrak{S}_n$ be a permutation. Suppose that $s$ is the unique positive integer satisfying $\omega(s)=n$. Define the \emph{trace} $\mathrm{Tr}_n (\omega)\in\mathfrak{S}_{n-1}$ of $\omega$ by 
    \[
        \mathrm{Tr}_n(\omega) (x) =
        \begin{cases}
            \omega (x), & x\neq s\\
            \omega (n), & x=s
        \end{cases}
    \]
    for every $1\leq x\leq n-1$. For simplicity, we will also write $\mathrm{Tr}:\mathfrak{S}_n\to\mathfrak{S}_{n-1}$ if $n$ is clear from the context.
\end{definition}
In other words, if the permutation $\omega$ is presented as a diagram from right to left, then the operation $\mathrm{Tr}:\mathfrak{S}_n\to\mathfrak{S}_{n-1}$ is exactly closing up the bottom two ends of $\omega$. 
\begin{example}
    Let $\omega=(243)\in\mathfrak{S}_4$ presented as 
    \[
        \omega = 
        \vcenter
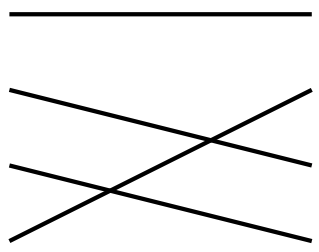
{\hbox{
            \begin{tikzpicture}[
                thick,
            ]
            \draw (0,0)--(2,0);
            \draw (0,-0.5)--(2,-1);
            \draw (0,-1)--(2,-1.5);
            \draw (0,-1.5)--(2,-0.5);             
            \end{tikzpicture}
        }}.
    \]
    Then $\mathrm{Tr}(\omega)=(23)\in\mathfrak{S}_3$ is presented as 
    \[
        \mathrm{Tr}(\omega) = 
        \vcenter
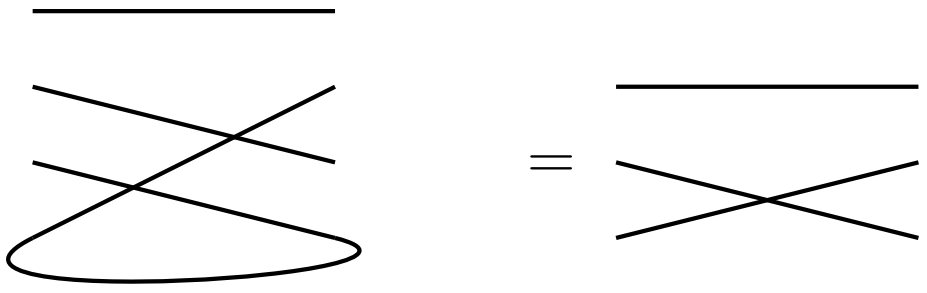
{\hbox{
            \begin{tikzpicture}[
                thick,
            ]
            \draw (0,0)--(2,0);
            \draw (0,-0.5)--(2,-1);
            \draw (0,-1)--(2,-1.5);
            \draw (0,-1.5)--(2,-0.5);    
            \draw (0,-1.5) ..controls +(-1,-0.5) and +(1,-0.25)..(2,-1.5);
            \end{tikzpicture}
        }}
        =
        \vcenter{\hbox{
            \begin{tikzpicture}[
                thick,
            ]
            \draw (0,0)--(2,0);
            \draw (0,-0.5)--(2,-1);
            \draw (0,-1)--(2,-0.5);             
            \end{tikzpicture}
        }}
        \in\mathfrak{S}_3.
    \]
\end{example}
\begin{remark}
    Write $e_k\in\mathfrak{S}_{k}$ as the identity permutation and $k$ is always omitted if it is clear from the context. Moreover, define 
    \[
        \begin{cases}
            A:=aq^{-2}\\
            Q:=q^2\\
            T:=t^2 q^{-2}\\
        \end{cases}
    \]
    where $a$, $q$, and $t$ are Hochschild grading, quantum grading and homological grading of the colored HOMFLY homology, respectively. The variables $(A,Q,T)$ are called the \emph{geometric variables}, while $(a,q,t)$ are called the \emph{algebraic variables}. The recursive formulae can be better expressed in geometric variables.
\end{remark}

Based on~\cite{conners2024row} and~\cite{conners2025doubly}, the recursive formulae for dimensions of reduced homology of $(\wedge^k,\wedge^1,\wedge^1,\ldots,\wedge^1)$-colored positive torus links are as follows:

\begin{theorem}[{\cite[Lemma 5.18]{conners2024row}}]\label{thm:recursive-formula-for-y-ified-homology}
For each pair $v \in {\{0, 1\}}^{l + m}$, $w \in {\{0, 1\}}^{l + n}$ and each $\sigma \in \mathfrak{S}^l$, let $p_{\sigma}(v, w) \in \mathbb{N}[[a, q, t]]$ be the family of power series satisfying the following properties:

\begin{enumerate}
    \item $p_e(0^m, \emptyset) = {\left( \cfrac{1 + A}{(1 - Q)(1 - T)} \right)}^m$ and $p_e(\emptyset, 0^n) = {\left( \cfrac{1 + A}{(1 - Q)(1 - T)} \right)}^n$;

    \item $p_e(0^m1, 0^n1) = \left( \cfrac{1 + A}{(1 - Q)(1 - T)} \right) p_e(0^m, 0^n)$ for all $m, n \geq 0$;

    \item $p_{\sigma}(v1, w1) = \left( \cfrac{Q^l + A}{1 - Q} \right) p_{Tr_{l + 1}(\sigma)}(v, w)$ if $|v| = |w| = l \geq 1$ and $\sigma(n) = n$;

    \item $p_{\sigma}(v1, w1) = \left( Q^l + A \right) p_{Tr_{l + 1}(\sigma)}(v, w)$ if $|v| = |w| = l \geq 1$ and $\sigma(n) \neq n$;

    \item $p_{\sigma}(v0, w1) = p_{\sigma {\pi_{(10^{l - 1})}}}(v, 1w)$ if $|v| = |w| + 1 = l \geq 1$;

    \item $p_{\sigma}(v1, w0) = p_{{\pi_{(10^{l - 1})}} \sigma}(1v, w)$ if $|v| + 1 = |w| = l \geq 1$;

    \item $p_e(0^m, 0^n) = p_e(10^{m - 1}, 10^{n - 1})$ for all $m, n \geq 1$;

    \item $p_{\sigma}(v0, w0) = Q^{-l}p_{{e_1 \sqcup \sigma}}(1v, 1w) + Q^{-l}Tp_{\sigma}(0v, 0w)$ if $|v| = |w| = l \geq 1$.
\end{enumerate}

These relations uniquely determine $p_{\sigma}(v, w)$ for all possible input data $(v, w, \sigma)$. Moreover, for each $m,n,k\geq 1$, with $d=\gcd(m,n)$, let 
\[
    f(m):=\prod_{k=1}^{m}\frac{Q^{\frac{k}{2}}-Q^{-\frac{k}{2}}}{Q^{\frac{1}{2}}-Q^{-\frac{1}{2}}}\cdot\prod_{k=1}^{m}\frac{1+AQ^{1-k}}{1-Q^k}.
\]
Then the positive torus link $T(m,n)$ satisfies
\begin{equation*}
    \dim\left(\overline{\mathrm{HHH}}_{\bigwedge^k,\bigwedge^1,\ldots,\bigwedge^1}(T(m,n))\right) = {(1-T)}^d{f(k)}^{-1} p_e(v_k,w_k).
\end{equation*}
\end{theorem}

In particular, for $\bigwedge^m$-coloring positive torus knot $T(2,2p+1)$, the bottom row of its column-colored HOMFLY homology satisfies the following recursive formulae. 
\begin{corollary}\label{cor:recursive-formula-for-bottom-row}
    For each pair $v \in {\{0, 1\}}^{l + m}$, $w \in {\{0, 1\}}^{l + n}$ and each $\sigma \in \mathfrak{S}^l$, let $p_{\sigma}(v, w) \in \mathbb{N}[[a, q, t]]$ be the family of power series as defined in Theorem~\ref{thm:recursive-formula-for-y-ified-homology}. Define $p'_e(v,w):=p_e(v,w)|_{A=0}\in\mathbb{N}[[q,t]]$. Then $p'_e(v,w)$ satisfies the following recursive formulae
\begin{enumerate}
    \item $p'_e(0^m, \emptyset) = {\left( \cfrac{1}{(1 - Q)(1 - T)} \right)}^m$ and $p'_e(\emptyset, 0^n) = {\left( \cfrac{1}{(1 - Q)(1 - T)} \right)}^n$;
    
    \item $p'_e(0^m1, 0^n1) = \left( \cfrac{1}{(1 - Q)(1 - T)} \right) p'_e(0^m, 0^n)$ for all $m, n \geq 0$;
    
    \item $p'_{\sigma}(v1, w1) = \left(\cfrac{Q^l}{1 - Q}\right) p'_{\mathrm{Tr}_{l + 1}(\sigma)}(v, w)$ if $|v| = |w| = l \geq 1$ and $\sigma(n) = n$;
    
    \item $p'_{\sigma}(v1, w1) = Q^l p'_{\mathrm{Tr}_{l + 1}(\sigma)}(v, w)$ if $|v| = |w| = l \geq 1$ and $\sigma(n) \neq n$;
    
    \item $p'_{\sigma}(v0, w1) = p'_{\sigma {\pi_{(10^{l - 1})}}}(v, 1w)$ if $|v| = |w| + 1 = l \geq 1$;
        
    \item $p'_{\sigma}(v1, w0) = p'_{{\pi_{(10^{l - 1})}} \sigma}(1v, w)$ if $|v| + 1 = |w| = l \geq 1$;
        
    \item $p'_e(0^m, 0^n) = p'_e(10^{m - 1}, 10^{n - 1})$ for all $m, n \geq 1$;
        
    \item $p'_{\sigma}(v0, w0) = Q^{-l}p'_{{e_1 \sqcup \sigma}}(1v, 1w) + Q^{-l}Tp'_{\sigma}(0v, 0w)$ if $|v| = |w| = l \geq 1$.
\end{enumerate}
Moreover, for integers $m\geq 1$ and $p\geq 0$, the positive torus knot $T(2,2p+1)$ satisfies
\[
    \dim\left(\overline{\mathrm{HHH}}_{\bigwedge^m}(T(2, 2p+1))\right) |_{A=0} = (1-T){f(m)}^{-1}p'_{e}(1^m 0^m, 1^m 0^{2pm}).
\] 
\end{corollary}

Another important result comes from orbit tracing. Let us first define a permutation associated to the tuple $(v,w,\sigma)$.

\begin{definition}
    For a positive integer $l$, let $v$ and $w$ be two $(01)$-sequences of length $l+m$ and $l+n$, where $v$ has $l$ zeros and $m$ ones, and $w$ has $l$ zeros and $n$ ones. Let $\sigma\in\mathfrak{S}^l$ be a permutation. Define the permutation associated to the tuple $(v,w,\sigma)$ as
    \begin{align*}
        {\omega(v,w)}_\sigma 
        &:= (e_n\sqcup\pi_v) \cdot (\pi_{1^m 0^n}\sqcup\sigma)\cdot (e_m\sqcup\pi^{-1}_w) \in \mathfrak{S}_{l+m+n}\\
        &=
        \vcenter
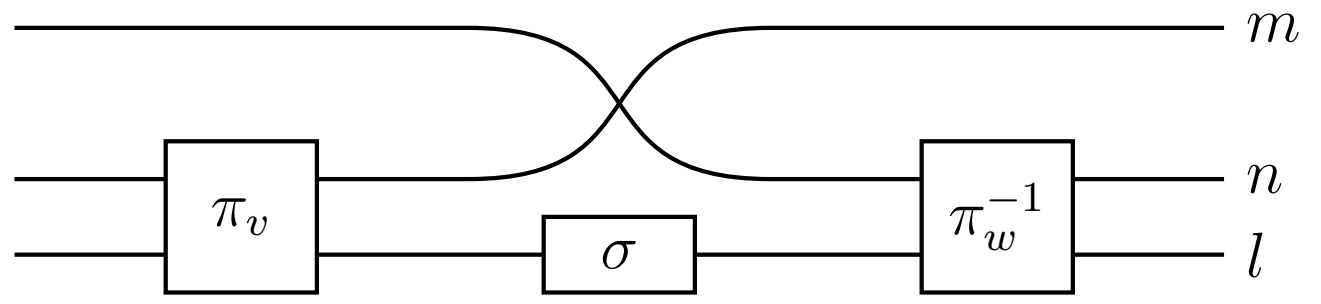
{\hbox{
            
\begin{tikzpicture}[
    baseline = (current bounding box.center),
]
    \begin{knot}[
        clip width = 5,
        looseness = 1.5,
        draft mode = strands,
    ]
        \strand[thick] (0,0.5) to [out = 0, in = 180] (8,0.5);
        \strand[thick] (0,1) to [out = 0, in = 180] (3,1)
        to [out = 0, in = 180] (5,2)
        to [out = 0, in = 180] (8,2);
        \strand[thick] (0,2) to [out = 0, in = 180] (3,2)
        to [out = 0, in = 180] (5,1)
        to [out = 0, in = 180] (8,1);
    \end{knot}

    \draw[thick, fill = white] (1,0.25) rectangle (2,1.25) node[midway] {$\pi_v$};
    \draw [thick, fill = white] (6,0.25) rectangle (7,1.25)
    node[midway] {$\pi^{-1}_w$};
    \draw [thick, fill = white] (3.5,0.25) rectangle (4.5,0.75)
    node[midway] {$\sigma$};

    \node[right] at (8,0.5) {$l$};
    \node[right] at (8,1) {$n$};
    \node[right] at (8,2) {$m$};
\end{tikzpicture}

        }},
    \end{align*}
    where $\pi_v$ and $\pi_w$ are the permutations associated to $v$ and $w$, respectively.
\end{definition}


In the following, we will focus on the number of orbits of the permutation ${\omega(v,w)}_\sigma$, which follows certain patterns under recursive formulae in Lemma~\ref{thm:recursive-formula-for-y-ified-homology}. Such pattern plays an essential role in tracking exponents in the formal power series.

\begin{lemma}[{\cite[Lemma 4.5]{conners2025doubly}}]\label{lemma:orbit-tracing}
    The permutation ${\omega(v,w)}_\sigma$ satisfies the following properties:
    \begin{enumerate}
        \item ${\omega(v_k,w_k)}_e$ has exactly $dk$ orbits for all $k$, $m$, and $n$, where $d=\gcd(m,n)$. In particular, ${\omega(1^m 0^m,1^m 0^{2pm})}_e$ has exactly $m$ orbits;
        \item Formulae (1), (2), and (3) in Theorem~\ref{thm:recursive-formula-for-y-ified-homology} decrease the number of orbits by $1$, while (4)-- (8) in Theorem~\ref{thm:recursive-formula-for-y-ified-homology} does not change the number of orbits. To be precise, 
        \begin{itemize}
            \item \(
            c({\omega(0^m,\emptyset)}_{e})=c({\omega(\emptyset,0^n)}_e)=1\);
            \item \(c({\omega(0^m 1,0^n 1)}_e)= c({\omega(0^m,0^n)}_e)+1\);
            \item \(c\left({\omega(v1,w1)}_\sigma\right)=c\left({\omega(v,w)}_{\mathrm{Tr}(\sigma)}\right)+1\), if \(\sigma(l)=l\), and $|v|+1=|w|+1=l$;
            \item \(c\left({\omega(v1,w1)}_\sigma\right)=c\left({\omega(v,w)}_{\mathrm{Tr}(\sigma)}\right)\), if \(\sigma(l)\neq l\), and $|v|+1=|w|+1=l$;
            \item \(c\left({\omega(v0,w1)}_\sigma\right)=c\left({\omega(v,1w)}_{\sigma\pi_{10^{l-1}}}\right)\), if $|v|=|w|+1=l\geq 1$;
            \item \(c\left({\omega(v1,w0)}_\sigma\right)=c\left({\omega(1v,w)}_{\pi_{10^{l-1}}\sigma}\right)\), if $|v|+1=|w|=l\geq 1$;
            \item \(c\left({\omega(0^m,0^n)}_e\right)=c\left({\omega(0^m 1,0^n 1)}_{e}\right)\);
            \item \(c\left({\omega(v0,w0)}_\sigma\right)
            =c\left({\omega(1v,1w)}_{e\sqcup\sigma}\right)
            =c\left({\omega(0v,0w)}_{\sigma}\right)\),
        \end{itemize}
        where $c(\omega)$ denotes the number of orbits in the permutation $\omega$;
    \end{enumerate}
\end{lemma}

\subsection{Elimination process}
From now on, we apply Corollary~\ref{cor:recursive-formula-for-bottom-row} to give an explicit formula for the bottom row of the column-colored homology, to prove Theorem~\ref{thm:main-theorem-explicit-formula}, and compare it with the conjecture in~\cite{gorsky2013quadruply}. Using Corollary~\ref{cor:recursive-formula-for-bottom-row}, the formal power series can be simplified through the following elimination process. 

To begin with, let us introduce the convention to apply formula (8) of Corollary~\ref{cor:recursive-formula-for-bottom-row}.

\begin{definition}[$\kappa$-admissible sequence]
    For fixed integers $m\geq 1$ and $p\geq 0$, consider a sequence of integers $\kappa:=(k_0,k_1,\ldots,k_p)$ with $m=k_0\geq k_1\geq k_2\geq\cdots\geq k_p\geq 0$. A sequence $(u_0,u_1,u_2,\ldots,u_p)$ is called \emph{$\kappa$-admissible}, if 
    \begin{enumerate}
        \item $u_0=0^m$, 
        \item $u_i\in {\{0,1\}}^{k_{i-1}}$, for every $1\leq i\leq p$, and 
        \item the number of zeros in $u_i$ is equal to $k_i$, for every $1\leq i\leq p$.
    \end{enumerate}
    Let $u_i(j)$ denote the $j$-th entry of $u_i$ for every $0\leq i\leq p$. Define the corresponding sequence ${\{u'_i\}}_{0\leq i\leq p}$ of $(u_0,u_1,\ldots,u_p)$ by the following recursive formula: 
    \begin{enumerate}
        \item $u'_0=u_0=0^m$, $u'_i\in{\{0,1\}}^m$, for every $1\leq i\leq p$.
        \item $u'_{i+1}$ is obtained from $u'_i$ by replacing the $j$-th zero in $u'_i$ by $u_{i+1}(j)$, for each $1\leq j\leq k_i$. 
    \end{enumerate}
\end{definition}

\begin{example}
    Let $(m,p)=(4,5)$, and $\kappa=(4,4,3,3,2,1)$. Then the following $\{u_i\}$ is an admissible $\kappa$-sequence:
    \begin{align*}
        u_0 &= 0000\\
        u_1 &= 0000\\
        u_2 &= 0100\\
        u_3 &= 000\\
        u_4 &= 010\\
        u_5 &= 10
    \end{align*}
    The corresponding $u'$-sequence is the following: 
    \begin{align*}
        u'_0 &= 0000,\quad k_0=4\\
        u'_1 &= 0000,\quad k_1=4\\
        u'_2 &= 0100,\quad k_2=3\\
        u'_3 &= 0100,\quad k_3=3\\
        u'_4 &= 0110,\quad k_4=2\\
        u'_5 &= 1110,\quad k_5=1
    \end{align*}
    Observe that the number of zeros in $u_i$ is equal to that in $u'_i$, where both are $k_i$.
\end{example}

With the $u$-sequence above, we are able to simplify the formal power series $p_e(1^m 0^m,\\1^m 0^{2pm})$ as the sum over all admissible sequences. The elimination process is described as follows: whenever both $v$ and $w$ end with $0$, rule (8) in Corollary~\ref{cor:recursive-formula-for-bottom-row} splits the formal power series into a sum of two terms:
\[
    p'_{\sigma}(v0, w0) = Q^{-l} p'_{e_1 \sqcup \sigma}(1v, 1w) + Q^{-l}T p'_{\sigma}(0v, 0w).
\]
In both terms on the right-hand side of the equation, the rightmost $0$ is moved to the leftmost position, where the first term corresponds to turning the rightmost $0$ into a $1$, and the second term corresponds to keeping the rightmost $0$ as a $0$. Each application of rule (8) thus encodes a binary choice at each such digits: when a trailing $0$ shows up in both $v$ and $w$, we decide whether it becomes a $1$ or remains a $0$ after being moved to the leftmost position. 

Observe that the initial stage of $(v,w)$ is $(1^m 0^m, 1^m 0^{2pm})$. Since $v$ and $w$ both have at least $m$ trailing zeros, then we start by making such a binary choice for $m$ times. In addition, we can encode the choices of $0$ and $1$ in a $u_1$-sequence, where $u_1$ is a $(01)$-sequence of length $m$ with $k_1$ zeros, and $u_1(i)$ is the choice for the $i$-th digit $0$. Using other rules in Corollary~\ref{cor:recursive-formula-for-bottom-row}, the stage can be reduced to 
\(
    (1^m u'_1 , u_1 1^m 0^{(2p-2)m})
\),
where $u'_1$ is obtained from $u_1$ by replacing the $j$-th zero in $u'_1$ by $u_1(j)$. In this new stage, we need to choose make $k_1$ choices, or equivalently a choice of $u_2$ of length $k_1$. We denote by $k_2$ the number of zeros in $u_2$. This procedure can be done for subsequent stages until the $p$-th stage
\(
    (1^m u'_p , u_p u_{p-1}\cdots u_1 1^m)
\),
where a choice $u_p$ of length $k_{p-1}$ with $k_p$ zeros is made, and the rightmost $2pm$ zeros in $w_m$ are eliminated. It turns out no further choices are needed in the rest of the computation. Hence, all the binary choices for the trailing $0$ can be encoded as a sequence $(u_1,u_2,\ldots,u_p)$, where each choice contributes to a term in the summation.

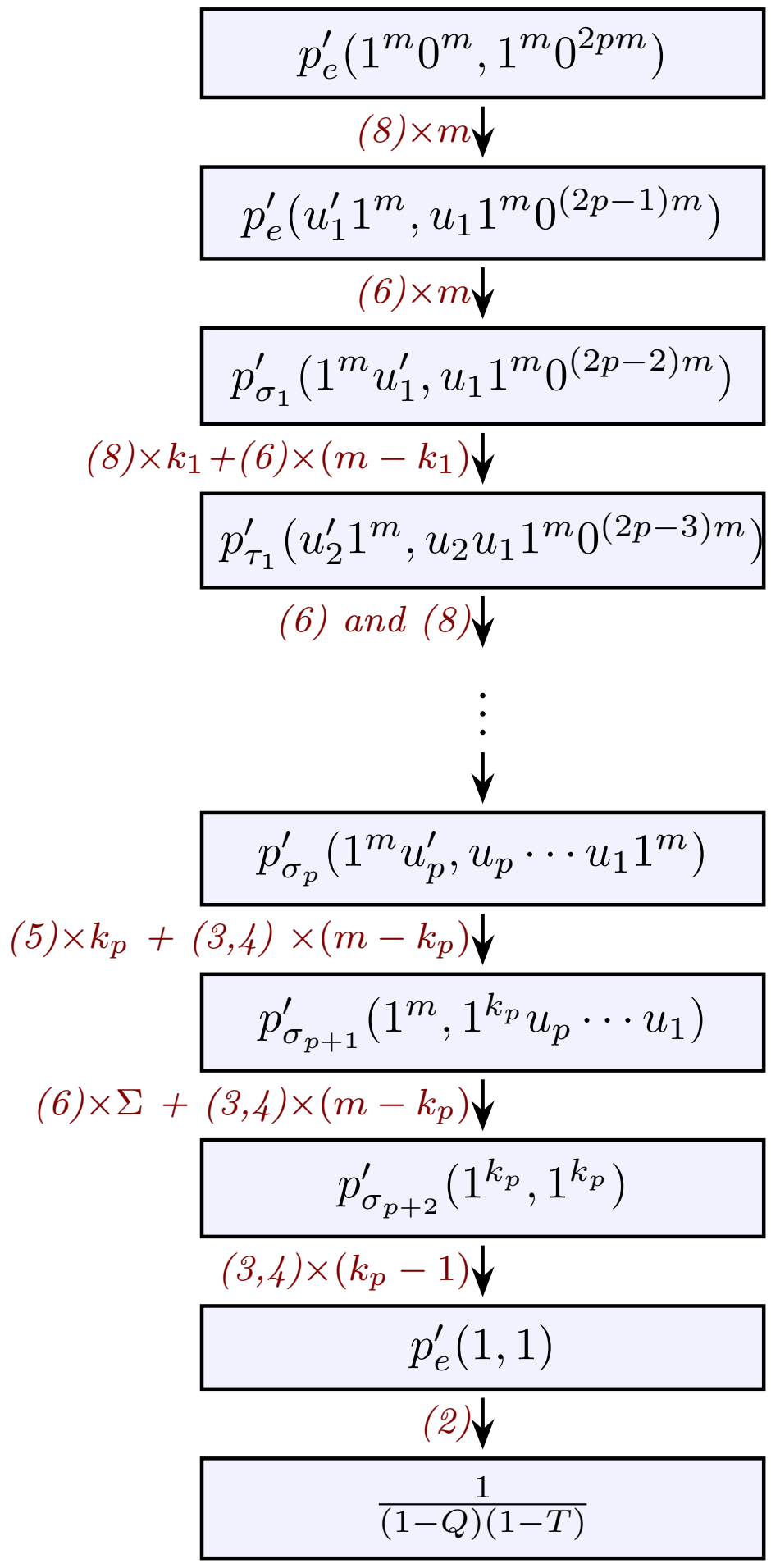
\begin{figure}[ht]
    \centering
    \begin{tikzpicture}[
  node distance = 0.5cm,
  poly/.style = {
    draw,
    thick,
    fill = blue!5,
    align = center,
    font = \footnotesize,
    inner sep = 4pt,
    text width = 4cm
  },
  arrow/.style = {
    ->,
    >=Stealth,
    thick,
    shorten >=1.5pt,
    shorten <=1.5pt,
  },
  method/.style = {
    text = red!50!black,
    font = \tiny\itshape,
    inner sep = 2pt,
  },
  dots/.style = {
    font = \normalsize,
    inner sep = 2pt,
  },
]

\node[poly] (P0) {
  $\displaystyle
    p'_e(1^m0^m, 1^m0^{2pm})
  $
};

\node[poly, below = of P0] (P1) {
  $\displaystyle
    p'_e(u_1'1^m,u_1 1^m 0^{(2p-1)m})
  $
};

\node[poly, below = of P1] (P2) {
  $\displaystyle
    p'_{\sigma_1}(1^mu_1',u_1 1^m 0^{(2p-2)m})
  $
};

\node[poly, below = of P2] (P3) {
  $\textstyle
    p'_{\tau_1}(u_2'1^m , u_2 u_1 1^m 0^{(2p-3)m})
  $
};

\node[below = of P3, dots] (D1) {$\vdots$};

\node[poly, below = of D1] (P4) {
  $\textstyle
  p'_{\sigma_p}(1^m u_{p}', u_p \cdots u_1 1^m)
  $
};

\node[poly, below = of P4] (P5) {
  $\textstyle
    p'_{\sigma_{p+1}}(1^m, 1^{k_p} u_p \cdots u_1)
  $
};

\node[poly, below = of P5] (D2) {
  $\textstyle
    p'_{\sigma_{p+2}} (1^{k_p},1^{k_p})
  $
};

\node[poly, below = of D2] (P6) {
  $\textstyle
    p'_e(1,1)
  $
};

\node[poly, below = of P6] (P7) {
  $\textstyle\frac{1}{(1-Q)(1-T)}$
};

\draw[arrow] (P0.south) -- (P1.north)
  node[midway, left, method] {(8)$\times m$};

\draw[arrow] (P1.south) -- (P2.north)
  node[midway, left, method] {(6)$\times m$};

\draw[arrow] (P2.south) -- (P3.north)
  node[midway, left, method] {(8)$\times k_1$+(6)$\times (m-k_1)$};

\draw[arrow] (P3.south) -- (D1.north)
  node[midway, left, method] {(6) and (8)};

\draw[arrow] (D1.south) -- (P4.north)
  node[midway, left, method] {};

\draw[arrow] (P4.south) -- (P5.north)
  node[midway, left, method] {(5)$\times k_p$ + (3,4) $\times (m-k_p)$};

\draw[arrow] (P5.south) -- (D2.north)
  node[midway, left, method] {(6)$\times\Sigma$ + (3,4)$\times (m-k_p)$};

\draw[arrow] (D2.south) -- (P6.north)
  node[midway, left, method] {(3,4)$\times (k_p-1)$};

\draw[arrow] (P6.south) -- (P7.north)
  node[midway, left, method] {(2)};

\end{tikzpicture}
    \caption{Elimination process}\label{fig:elimination-process}
\end{figure}

Formally, the process can be written as follows for $u_1=(u_1(1)u_1(2)\cdots u_1(m)) \in {\{0,1\}}^m$:
\begin{align*}
    p'_e(1^m 0^m, 1^m 0^{2pm}) 
    =& Q^{-m} p'_e(1~1^m 0^{m-1},1~1^m 0^{2pm-1}) + Q^{-m} T p'_e(0~1^m0^{m-1}, 0~1^m 0^{2pm-1})\\
    =& \sum_{u_1(m)\in {\{0,1\}}} Q^{-m} T^{\mathbf{1}_{u_1(m)=0}} p'_e\big(u'_1(m)1^m 0^{m-1},u_1(m)1^m 0^{2pm-1}\big)\\
    =& \sum_{u_1(m-1)\in {\{0,1\}}} \sum_{u_1(m)\in {\{0,1\}}} Q^{-2m-\mathbf{1}_{u_1(m)=1}} T^{\mathbf{1}_{u_1(m-1)=0}+\mathbf{1}_{u_1(m)=0}}\times\\
    &p'_e\big(u'_1(m-1)u'_1(m)1^m 0^{m-2}, u_1(m-1)u_1(m)1^m 0^{2pm-2}\big)\\
    =&\cdots=\sum_{u_1\in {\{0,1\}}^m} Q^{\alpha_1} T^{k_1} p'_e(u'_1 1^m, u_1 1^m 0^{2pm-m}),
\end{align*}
where $\alpha_1$ is a constant depending on the choice of $u_1$, and $k_1$ is the number of zeros in $u_1$. In general, using (8) and (6), 
\begin{align}
    p'_e(1^m 0^m, 1^m 0^{2pm})
    &=\sum_{u_1\in {\{0,1\}}^m} Q^{\alpha_1} T^{k_1} p'_e(u'_1 1^m, u_1 1^m 0^{(2p-1)m})\nonumber\\
    &=\sum_{u_1\in {\{0,1\}}^m} Q^{\alpha_1} T^{k_1} p'_{\sigma_1} (1^m u'_1, u_1 1^m 0^{(2p-2)m})\nonumber\\
    &=\sum_{u_2\in{\{0,1\}}^{k_1}} \sum_{u_1\in {\{0,1\}}^m} Q^{\alpha_2} T^{k_1 +k_2} p'_{\tau_1}(u'_2 1^m, u_2 u_1 1^m 0^{(2p-3)m}) \nonumber\\
    &=\sum_{u_2\in{\{0,1\}}^{k_1}} \sum_{u_1\in {\{0,1\}}^m} Q^{\alpha_2} T^{k_1 +k_2} p'_{\sigma_2}(1^m u'_2, u_2 u_1 1^m 0^{(2p-4)m}) \nonumber\\
    &=\cdots=\mathop{\sum\cdots\sum}_{u_i\in{\{0,1\}}^{k_{i-1}},1\leq i\leq p} Q^{\alpha_p} T^{k_1+\cdots+k_p} 
    p'_{\sigma_p}(1^m u'_p,u_p u_{p-1}\cdots u_1 1^m),
\end{align}
where $\alpha_i$'s are constants depending on $u_i$, and $\sigma_i$'s are permutations depending on $u_i$. Through this process, we eliminate the $0^{2pm}$ in the second entry. Using (3--5) and (6), we have 
\begin{align}
    p'_e(1^m 0^m, 1^m 0^{2pm}) 
    &=\mathop{\sum\cdots\sum}_{u_i\in {\{0,1\}}^{k_{i-1}},1\leq i\leq p} \frac{Q^{\alpha_{p+1}} T^{k_1+\cdots+k_p}}{{(1-Q)}^{\gamma_1}} 
    p'_{\sigma_{p+1}}(1^m,1^{k_p}u_p u_{p-1}\cdots u_1)\nonumber\\
    &=\mathop{\sum\cdots\sum}_{u_i\in {\{0,1\}}^{k_{i-1}},1\leq i\leq p} \frac{Q^{\alpha_{p+2}} T^{k_1+\cdots+k_p}}{{(1-Q)}^{\gamma_2}}
    p'_{\sigma_{p+2}}(1^{k_p},1^{k_p})\nonumber\\
    &=\mathop{\sum\cdots\sum}_{u_i\in {\{0,1\}}^{k_{i-1}},1\leq i\leq p} \frac{Q^{\alpha} T^{k_1+\cdots+k_p}}{{(1-Q)}^{\gamma_3}}
    p'_e(1,1)\nonumber\\
    &\xlongequal{\text{Rule (2)}}\mathop{\sum\cdots\sum}_{u_i\in{\{0,1\}}^{k_{i-1}},1\leq i\leq p} \frac{Q^{\alpha} T^{k_1+\cdots+k_p}}{{(1-Q)}^{\gamma_3+1}(1-T)},\label{eq:pe1}
\end{align}
where $\gamma_1$, $\gamma_2$, and $\gamma_3$ are constants depending on the choices of $u_i$'s. The whole process can be depicted as Figure~\ref{fig:elimination-process}.

The following example illustrates how to obtain a term in the summation of Equation~\ref{eq:pe1} through the elimination process.
\begin{example}
    Let $(m,p)=(2,2)$. Then $p'_e(1^m 0^m,1^m 0^{2pm})=p'_e(1100,1100000000)$. Denote by ${(v,w)}_\sigma$ the power series $p'_\sigma(v,w)$. By the elimination process described above, we need to first make $m=2$ binary choices to the two trailing zeros in $v$ and $w$. Choose $u_1:=(10)$. Applying rule (8) in Corollary~\ref{cor:recursive-formula-for-bottom-row}, we have
    \[
        {(11\underline{00},11 00000000)}_{e_2} \to {(10 11, 10 11 000000)}_{e_3},
    \]
    where the binary choices are maded to the trailing zeros in the underlined digits of $v$. In this case, $k_1=1$ and $u'_1=(10)$ as $u_1$ has exactly one zero. Applying rule (6) in Corollary~\ref{cor:recursive-formula-for-bottom-row}, we obtain the first stage 
    \[
        {(1011,1011000000)}_{e_3} \overset{(6)}{\to} {(1110,10110000)}_{(132)}
    \]
    as $\pi_{100}\cdot\pi_{100}\cdot e_3=(132)$. In the second elimination, we need to make $k_1=1$ binary choice to the trailing zeros in $v$ and $w$. Choose $u_2:=(0)$. Applying rule (8) for the zero in the underlined term, and rule (6) for the ones in the underlined term, we have
    \[
        {(11\underline{10},10110000)}_{(132)} \to {(0111,01011000)}_{(132)}\overset{(6)}{\to}{(1011,0101100)}_{e_3},
    \]
    as $\pi_{100}\cdot (132)=e_3$. Then $k_2=1$, and $u'_2=(10)$. Apply rule (6) twice to the above term, we obtain the second stage
    \[
        {(1011,0101100)}_{e_3} \overset{(6)}{\to} {(1110,01011)}_{(132)} = {(11 u'_2, u_2 u_1 11)}_{(132)},
    \]
    which is the $p$-th stage of the elimination process. Since no further binary choices are required, then $\kappa=(k_0,k_1,k_2)=(2,1,1)$, and the above term is the summand associated to the $\kappa$-admissible sequence $(u_0,u_1,u_2)=(00,10,0)$ in the summation of Equation~\ref{eq:pe1}. 
\end{example}
In conclusion, the elimination process reduces the computation of $p_e(1^m 0^m, 1^m 0^{2pm})$ to a sum over all admissible sequences $(u_0,u_1,u_2,\ldots,u_p)$, according to Equation~\ref{eq:pe1}. In the following sections, we will compute the exponents $\alpha$ and $\gamma_3$ in Equation~\ref{eq:pe1}, and give an explicit formula for the bottom row of the column-colored HOMFLY homology of $T(2,2p+1)$.
\subsection{Exponent of \texorpdfstring{$(1-Q)$}{(1-Q)}}
    By Lemma~\ref{lemma:orbit-tracing} (2), the number of orbits in $\omega(1^m 0^m,1^m 0^{2pm})$ is $m$. Note that in the elimination process, only (2), (3), (4), (5), (6), and (8) in Corollary~\ref{cor:recursive-formula-for-bottom-row} are applied, among which only (2) and (3) change the number of orbits. In addition, the number of orbits decreases by $1$, exactly when a (2) or (3) is applied. Rule (2) is applied exactly once for every $(u_0,u_1,u_2,\ldots,u_p)$. Therefore, the exponent of $(1-Q)$ is equal to minus the number of times that rule (2) and (3) are applied in the process, which is equal to the number of orbit changes through the process. 

    Since the final stage ${\omega(1,1)}_e$ contains precisely one orbit, and the initial stage ${\omega(1^m 0^m,\newline 1^m 0^{2pm})}_e$ contains $m$ orbits, then the net change in the number of orbits arising from (3) equals  $(m-1)$, while the contribution of (2) is $1$. Thus, the exponent of $(1-Q)$ is equal to $-m$. In conclusion, 
    \begin{equation}
        p'_e(1^m 0^m, 1^m 0^{2pm}) = {(1-Q)}^{-m}\mathop{\sum\sum\cdots\sum}_{0\leq k_p\leq k_{p-1}\leq\cdots\leq k_1\leq m} 
        \frac{T^{k_1+\cdots+k_p}}{1-T}\left(\sum_{\#_0 u'_j=k_j} Q^{\alpha} \right).\label{eq:pe2}
    \end{equation}
    
\subsection{Exponent of \texorpdfstring{$Q$}{Q}}
In this section, we compute the exponent of $Q$ in the summand of equation. 
Notice that during the simplification, only rule (3), (4), and (8) will change the exponent of $Q$ by the number of ones in $v$ or $w$, in which the number of ones are identical. Notice that (3) and (4) provides the same number of exponent of $Q$. Thus, it suffices to trace the number of $1$'s in the elimination. 

The computation splits into two parts, namely $\alpha=E_1+E_2$, where $E_1$ denotes the exponent changes of $Q$ when transforming from $p_e(1^m 0^m, 1^m 0^{2pm})$ to $p_{\sigma_p}(1^m u'_p\cdots u'_1 1^m)$, and $E_2$ denotes the exponent changes of $Q$ when transforming from $p_{\sigma_p}(1^m u'_p\cdots u'_1 1^m)$ to $p_e(1,1)$. 

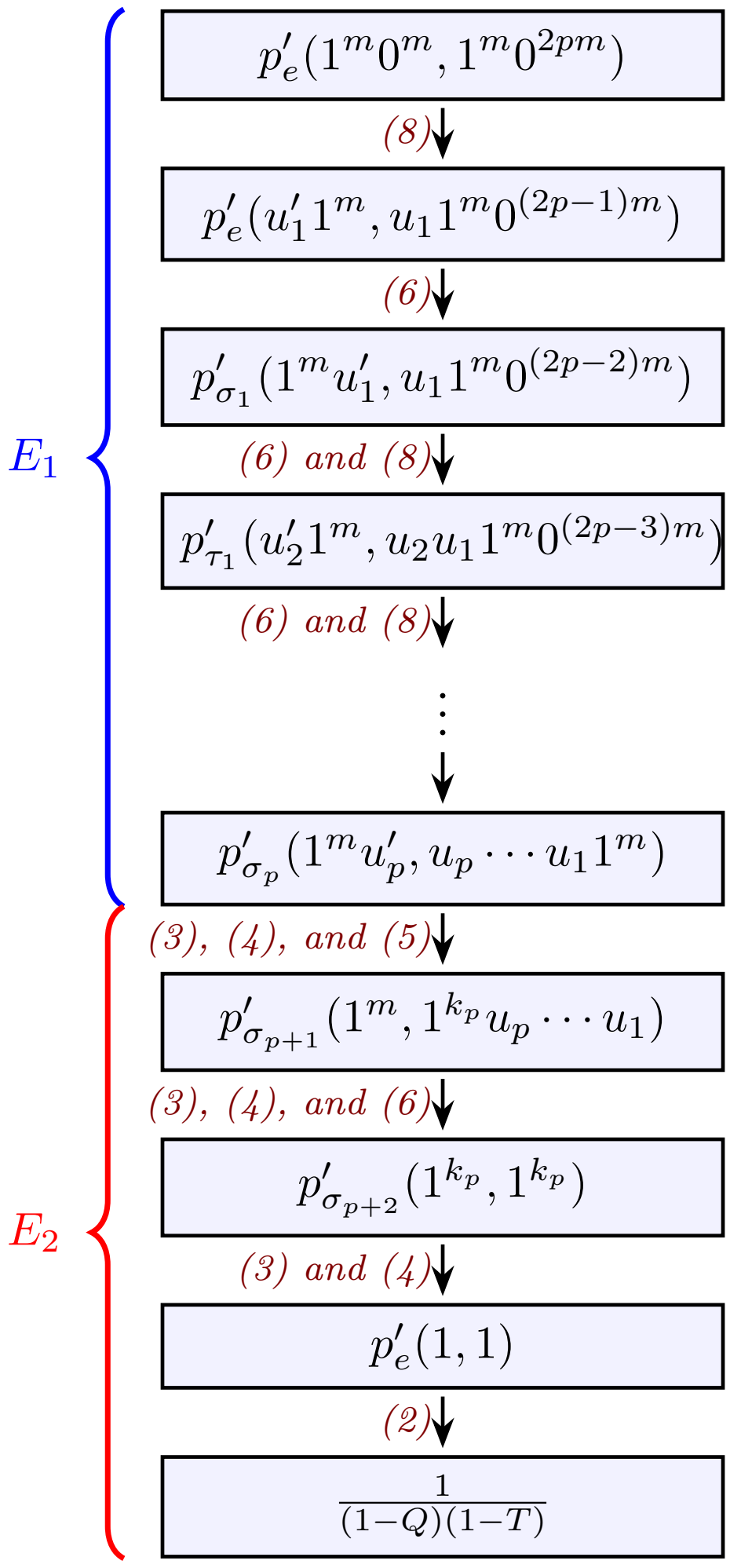
\begin{figure}[ht]
    \centering
    \begin{tikzpicture}[
  node distance = 0.5cm,
  poly/.style = {
    draw,
    thick,
    fill = blue!5,
    align = center,
    font = \footnotesize,
    inner sep = 4pt,
    text width = 4cm
  },
  arrow/.style = {
    ->,
    >=Stealth,
    thick,
    shorten >=1.5pt,
    shorten <=1.5pt,
  },
  method/.style = {
    text = red!50!black,
    font = \tiny\itshape,
    inner sep = 2pt,
  },
  dots/.style = {
    font = \normalsize,
    inner sep = 2pt,
  },
]

\node[poly] (P0) {
  $\displaystyle
    p'_e(1^m0^m, 1^m0^{2pm})
  $
};

\node[poly, below = of P0] (P1) {
  $\displaystyle
    p'_e(u_1'1^m,u_1 1^m 0^{(2p-1)m})
  $
};

\node[poly, below = of P1] (P2) {
  $\displaystyle
    p'_{\sigma_1}(1^mu_1',u_1 1^m 0^{(2p-2)m})
  $
};

\node[poly, below = of P2] (P3) {
  $\textstyle
    p'_{\tau_1}(u_2'1^m , u_2 u_1 1^m 0^{(2p-3)m})
  $
};

\node[below = of P3, dots] (D1) {$\vdots$};

\node[poly, below = of D1] (P4) {
  $\textstyle
  p'_{\sigma_p}(1^m u_{p}', u_p \cdots u_1 1^m)
  $
};

\node[poly, below = of P4] (P5) {
  $\textstyle
    p'_{\sigma_{p+1}}(1^m, 1^{k_p} u_p \cdots u_1)
  $
};

\node[poly, below = of P5] (D2) {
  $\textstyle
    p'_{\sigma_{p+2}} (1^{k_p},1^{k_p})
  $
};

\node[poly, below = of D2] (P6) {
  $\textstyle
    p'_e(1,1)
  $
};

\node[poly, below = of P6] (P7) {
  $\textstyle\frac{1}{(1-Q)(1-T)}$
};

\draw[arrow] (P0.south) -- (P1.north)
  node[midway, left, method] {(8)};

\draw[arrow] (P1.south) -- (P2.north)
  node[midway, left, method] {(6)};

\draw[arrow] (P2.south) -- (P3.north)
  node[midway, left, method] {(6) and (8)};

\draw[arrow] (P3.south) -- (D1.north)
  node[midway, left, method] {(6) and (8)};

\draw[arrow] (D1.south) -- (P4.north)
  node[midway, left, method] {};

\draw[arrow] (P4.south) -- (P5.north)
  node[midway, left, method] {(3), (4), and (5)};

\draw[arrow] (P5.south) -- (D2.north)
  node[midway, left, method] {(3), (4), and (6)};

\draw[arrow] (D2.south) -- (P6.north)
  node[midway, left, method] {(3) and (4)};

\draw[arrow] (P6.south) -- (P7.north)
  node[midway, left, method] {(2)};

\draw [decorate, decoration={brace, amplitude=7pt, mirror}, very thick, blue]
              ([xshift=-8pt]P0.north west) -- ([xshift=-8pt]P4.south west)
              node [midway, left=9pt, font=\footnotesize\bfseries, blue] {$E_1$};
\draw [decorate, decoration={brace, amplitude=7pt, mirror}, very thick, red]
              ([xshift=-8pt]P4.south west) -- ([xshift=-8pt]P7.south west)
              node [midway, left=9pt, font=\footnotesize\bfseries, red] {$E_2$};
\end{tikzpicture}
    \caption{Exponent of $Q$}\label{fig:splitting-of-q-exponent}
\end{figure}
\FloatBarrier{}
    
    \subsubsection{Computation of \texorpdfstring{$E_2$}{E2}}
    The computation of $E_2$ is straightforward. Observe that from $p_{\sigma_p}(1^m u_p',u_p \cdots u_1 1^m)$ to $p_e(1,1)$, the number of ones in $v$ or $w$ reduces exactly when we apply (3) or (5), while no (8) is applied. Therefore, the exponent of $Q$ in the summand increases by $l$ if and only if the number of ones decreases by $1$, where $l$ is the number of ones in the resulting sequence $v$ (or $w$). Hence,
    \begin{align}
        E_2:
        &=\sum_{i=1}^{2m-k_p-1} i\nonumber\\
        &=\frac{(2m-k_p-1)(2m-k_p)}{2}\nonumber\\
        &=2m^2-(2k_p+1)m+\frac{1}{2}k_p^2+\frac{1}{2}k_p.\label{e2result}
    \end{align}

    \subsubsection{Computation of \texorpdfstring{$E_1$}{E1}}
    Now we consider the contribution of rule (8) when $p_e(1^m 0^m,\allowbreak 1^m 0^{2pm})$ is transformed to $p_{\sigma_p}(1^m u_p',u_p \cdots u_1 1^m)$. For the $(j+1)$-st round we apply rule (8), at the $k$-th position, the $v$-entry becomes 
    \[
        \underline{u'_{j+1}(k+1)u'_{j+1}(k+2)\cdots u'_{j+1}(m)} 1^m \underline{u'_j(1)u'_j(2)\cdots u'_j(k-1) u'_j(k)}
    \]
    If $u'_j(k)=0$, then we apply rule (8) and the contribution to the exponent of $Q$ is equal to 
    \[
        -\left(m+\sum_{i=1}^{k-1} \mathbf{1}_{u'_{j}(i)=1} + \sum_{i=k+1}^{m} \mathbf{1}_{u'_{j+1}(i)=1} \right),
    \]
    no matter which part we choose. If $u'_j(k)=1$, then we apply rule (6) with no contribution to the exponent of $Q$. In conclusion, in the whole process transforming $p_e(1^m 0^m, 1^m 0^{2pm})$ to $p_{\sigma_p}(1^m u_p',u_p \cdots u_1 1^m)$, the total contribution is 
    \[
        \begin{aligned}
            E_1:
            &=\sum_{j=0}^{p-1}\sum_{k:u'_j(k)=0} -\left(m+\sum_{i=1}^{k-1} \mathbf{1}_{u'_{j}(i)=1} + \sum_{i=k+1}^{m} \mathbf{1}_{u'_{j+1}(i)=1} \right)\\
            &=-\sum_{j=0}^{p-1} \left(m k_j+\sum_{k:u'_j(k)=0} \left(\sum_{i=1}^{k-1} \mathbf{1}_{u'_{j}(i)=1} + \sum_{i=k+1}^{m} \mathbf{1}_{u'_{j+1}(i)=1} \right)\right)\\
            &=-m^2-m\sum_{j=1}^{p-1} k_j - \sum_{j=0}^{p-1}\sum_{k:u'_j(k)=0} \left(\sum_{i=1}^{k-1} \mathbf{1}_{u'_{j}(i)=1} + \sum_{i=k+1}^{m} \mathbf{1}_{u'_{j+1}(i)=1} \right)
        \end{aligned}
    \]
    For a family of $(01)$-sequence $u'_0,u'_1,\ldots, u'_p$ of length $m$, where $u'_0=0^m$ and $u'_{i+1}$ is obtained from $u'_i$ by transforming some zeros in $u'_i$ to ones, and $u'_j$ has exactly $k_j$ zeros, for every $0\leq j\leq p$, we define its corresponding $\Delta$-sequence by
    \[
        \Delta(k):=\max_{0\leq j\leq p} \{j:u'_j(k)=0\}\in\{0,1,2,\ldots,p\}.
    \]
    Then $u'_j(i)=1$ if and only if $\Delta(i)\leq j-1$. Therefore, 
    \[
        \begin{aligned}
            E_1
            &=-m^2-m\sum_{j=1}^{p-1} k_j - \sum_{j=0}^{p-1}\sum_{k:u'_j(k)=0} \left(\sum_{i=1}^{k-1} \mathbf{1}_{\Delta(i)\leq j-1} + \sum_{i=k+1}^{m} \mathbf{1}_{\Delta(i)\leq j} \right)\\
            &= -m^2-m\sum_{j=1}^{p-1} k_j - 
            \sum_{k=1}^{m}\sum_{j=0}^{\Delta(k)} \mathbf{1}_{j\leq p-1}\left(\sum_{i=1}^{k-1}\mathbf{1}_{\Delta(i)\leq j-1} + \sum_{i=k+1}^{m}\mathbf{1}_{\Delta(i)\leq j}\right)\\
            &= -m^2-m\sum_{j=1}^{p-1} k_j - \sum_{k=1}^{m} \sum_{j=0}^{\Delta(k)} \mathbf{1}_{j\leq p-1} 
            \left(\sum_{i\neq k}\mathbf{1}_{\Delta(i)\leq j-1} + \sum_{i=k+1}^{m} \mathbf{1}_{\Delta(i)=j}\right).\\
        \end{aligned}
    \]
    Observe that $j-1<j\leq \Delta(k)$. Then $\mathbf{1}_{\Delta(k)\leq j-1}=0$, for every $0\leq j\leq \Delta(k)$. Therefore, 
        \begin{align}
            E_1
            &=-m^2-m\sum_{j=1}^{p-1} k_j - \sum_{k=1}^{m} \sum_{j=0}^{\Delta(k)} \mathbf{1}_{j\leq p-1} 
            \left(\sum_{i=1}^m\mathbf{1}_{\Delta(i)\leq j-1} + \sum_{i=k+1}^{m} \mathbf{1}_{\Delta(i)=j}\right)\nonumber\\
            &= -m^2-m\sum_{j=1}^{p-1} k_j - S_1 -S_2,\label{e1expression}
        \end{align}
    where 
    \[
        S_1:= \sum_{k=1}^{m} \sum_{j=0}^{\Delta(k)} \sum_{i=1}^m \mathbf{1}_{j\leq p-1} \mathbf{1}_{\Delta(i)\leq j-1},
    \]
    and 
    \[S_2:= \sum_{k=1}^{m} \sum_{j=0}^{\Delta(k)} \sum_{i=k+1}^m \mathbf{1}_{j\leq p-1}\mathbf{1}_{\Delta(i)=j}.\]
    For fixed $j$ and $k$, notice that $\Delta(i)\leq j-1$ if and only if $u'_j(i)=1$. Then 
    \[
        \sum_{i=1}^{m} \mathbf{1}_{\Delta(i)\leq j-1}=\#\{1\leq i\leq m: u'_j(i)=1\}=m-k_j,\quad \forall 1\leq j\leq p.
    \] 
    Therefore,
        \begin{align}
            S_1 
            &= \sum_{l=0}^{p} \sum_{k=1}^{m} \sum_{j=0}^{p-1} \mathbf{1}_{\Delta(k)\geq j} \mathbf{1}_{\Delta(k)=l} (m-k_j)\nonumber\\
            &=\sum_{l=0}^{p} \sum_{k=1}^{m} \left(\mathbf{1}_{\Delta(k)=l} \sum_{j=0}^{l} (m-k_j)\cdot\mathbf{1}_{j\leq p-1}\right)\nonumber\\
            &=\sum_{l=0}^{p-1} \sum_{k=1}^{m} \left(\mathbf{1}_{\Delta(k)=l} \left(ml-\sum_{j=1}^{l} k_j \right)\right)+ \sum_{k=1}^{m}\left(\mathbf{1}_{\Delta(k)=p}\sum_{j=0}^{p-1}(m-k_j)\right)\nonumber\\
            &= \sum_{l=0}^{p-1} \left(ml-\sum_{j=1}^{l} k_j \right) \sum_{k=1}^{m} \mathbf{1}_{\Delta(k)=l} 
            + \sum_{k=1}^{m} \mathbf{1}_{\Delta(k)=p}\cdot \left(m(p-1)-\sum_{j=1}^{p-1}k_j\right)\nonumber\\
            &= \sum_{l=0}^{p-1} (k_l-k_{l+1})\left(ml-\sum_{j=1}^{l} k_j \right)
            + k_p\cdot \left(m(p-1)-\sum_{j=1}^{p-1}k_j\right)\nonumber\\
            &= m\left(\sum_{l=0}^{p-1}lk_l - \sum_{l=0}^{p-1}l k_{l+1}\right)
            -\sum_{l=0}^{p-1} \left(k_l\sum_{j=1}^{l}k_j-k_{l+1}\sum_{j=1}^{l}k_j\right)
            +m(p-1)k_p-k_p\sum_{j=1}^{p-1} k_j\nonumber\\
            &=-mpk_{p} + m \sum_{l=0}^{p-1}k_{l+1} 
            + k_p\sum_{j=1}^{p} k_j - \sum_{l=0}^{p-1} k_{l+1}^2
            + mpk_p-mk_p-k_p\sum_{j=1}^{p-1} k_j\nonumber\\
            &= m(k_1+k_2+\cdots+k_{p-1}) -(k_1^2+k_2^2+\cdots +k_{p-1}^2).\label{s1result}
        \end{align}
    In addition, 
        \begin{align}
            S_2 
            &=\sum_{l=0}^{p} \sum_{k=1}^{m} \sum_{j=0}^{p-1} \mathbf{1}_{\Delta(k)=l} \cdot \mathbf{1}_{\Delta(k)\geq j} \sum_{i=k+1}^{m} \mathbf{1}_{\Delta(i)=j}\nonumber\\
            &=\sum_{l=0}^{p} \mathop{\sum\sum}_{1\leq k<i\leq m} \mathbf{1}_{\Delta(i)\leq p-1}\mathbf{1}_{\Delta(k)=l} \mathbf{1}_{\Delta(k)\geq\Delta(i)}\nonumber\\
            &=\mathop{\sum\sum}_{1\leq k<i\leq m} \mathbf{1}_{\Delta(i)\leq p-1}\mathbf{1}_{\Delta(k)\geq\Delta(i)}\nonumber\\
            &=\mathop{\sum\sum}_{1\leq k<i\leq m} \mathbf{1}_{\Delta(k)\geq\Delta(i)} - \sum_{\Delta(i)=p}\sum_{1\leq k< i} \mathbf{1}_{\Delta(i)=p}\nonumber\\
            &=\mathrm{Inv}(\Delta) -\binom{k_p}{2}.\label{s2result}
        \end{align}
    Combining~(\ref{e2result}),~(\ref{e1expression}),~(\ref{s1result}), and~(\ref{s2result}), we obtain
    \begin{align}
        E_1
        &= - \left(m^2 + m\sum_{j=1}^{p-1}k_j +m\sum_{j=1}^{p-1}k_j -\sum_{j=1}^{p-1}k_j^2+ \mathrm{Inv}(\Delta) -\binom{k_p}{2}\right)\nonumber\\
        &= -\left(m^2 + 2m\sum_{j=1}^{p-1}k_j -\sum_{j=1}^{p-1}k_j^2+ \mathrm{Inv}(\Delta) -\binom{k_p}{2}\right)\nonumber
    \end{align}
    and 
    \begin{align}
        \alpha
        &=E_1+E_2 \nonumber\\
        &=2m^2-(2k_p+1)m+\frac{1}{2}k_p^2+\frac{1}{2}k_p
        -\left(m^2 + 2m\sum_{j=1}^{p-1}k_j -\sum_{j=1}^{p-1}k_j^2+ \mathrm{Inv}(\Delta) -\binom{k_p}{2}\right)\nonumber\\
        &=m^2 - m - 2m\sum_{j=1}^{p}k_j + \sum_{j=1}^{p} k_j ^2 - \mathrm{Inv}(\Delta)\nonumber\\
        &=m^2 - m - 2m\sum_{j=1}^{p}k_j + \sum_{j=1}^{p} k_j ^2 - \mathrm{inv}(\Delta)-\sum_{j=0}^{p}\binom{k_j-k_{j+1}}{2}\nonumber\\
        &=m^2 - m - 2m\sum_{j=1}^{p}k_j + \sum_{j=1}^{p} k_j ^2 - \mathrm{inv}(\Delta)-\frac{1}{2}\sum_{j=0}^{p}(k_j^2+k_{j+1}^2-2k_j k_{j+1}-k_j+k_{j+1})\nonumber\\
        &=\frac{1}{2}m^2 -\frac{1}{2}m + \sum_{j=1}^{p}k_{j-1}k_j-2m\sum_{j=1}^{p}k_j-\mathrm{inv}(\Delta).\label{eq:alpha-result}
    \end{align}

\subsection{Final results}
    In conclusion, the explicit formula for the bottom row of the reduced colored HOMFLY homology of the $\bigwedge^m$-colored two-strand positive torus knot $T(2,2p+1)$ can be given as follows, which verifies the conjecture given by~\cite{gorsky2013quadruply}.

    \begin{theorem}\label{thm:main-theorem-explicit-formula}
        Let $m\geq 1$, $p\geq 0$ be two integers. Then the bottom row of reduced column-colored HOMFLY homology of $\bigwedge^m$-coloring positive torus knot $T(2,2p+1)$ is given by 
        \[
            \dim \left( \overline{\mathrm{HHH}}_{\bigwedge^m} (T(2,2p+1)) \right) |_{A=0} = 
            Q^{\frac{3m^2-3m}{4}}
            \mathop{\sum\sum\cdots\sum}_{0\leq k_p\leq k_{p-1}\leq\cdots\leq k_1\leq k_0=m}
            \sqbinom{m}{k_1}_{Q^{-1}}\sqbinom{k_1}{k_2}_{Q^{-1}}\times\cdots
        \]
        \[
            \times\sqbinom{k_{p-1}}{k_p}_{Q^{-1}} Q^{\sum\limits_{j=1}^{p}k_{j-1}k_j-2m\sum\limits_{j=1}^{p}k_j}T^{\sum\limits_{j=1}^{p}k_j}.
        \]
    \end{theorem}

    \begin{proof}
        Combining (\ref{eq:alpha-result}), and (\ref{eq:pe2}), using Lemma~\ref{lem:multinomial-theorem}, we have
        \begin{align}
            p'_e(1^m 0^m,1^m 0^{2pm})
            &= \mathop{\sum\sum\cdots\sum}_{0\leq k_p\leq k_{p-1}\leq\cdots\leq k_1\leq m} 
            \frac{T^{\sum\limits_{j=1}^{p}k_j}}{1-T}\left(\sum_{\Delta\in S_\kappa} \frac{Q^{\frac{1}{2}m^2 -\frac{1}{2}m + \sum\limits_{j=1}^{p}k_{j-1}k_j-2m\sum\limits_{j=1}^{p}k_j-\mathrm{inv}(\Delta)} }{{(1-Q)}^{m}}\right)\nonumber\\
            &=\frac{Q^{\frac{m^2-m}{2}}}{(1-T){(1-Q)}^m}
            \mathop{\sum\sum\cdots\sum}_{0\leq k_p\leq k_{p-1}\leq\cdots\leq k_1\leq m} 
            Q^{\sum\limits_{j=1}^{p}k_{j-1}k_j-2m\sum\limits_{j=1}^{p}k_j}T^{\sum\limits_{j=1}^{p}k_j}\nonumber\\
            &\hspace*{8em}\times\left(\sum_{\Delta\in S_\kappa} Q^{- \mathrm{inv}(\Delta)}\right)\nonumber\\
            &=\frac{Q^{\frac{m^2-m}{2}}}{(1-T){(1-Q)}^m}
            \mathop{\sum\sum\cdots\sum}_{0\leq k_p\leq k_{p-1}\leq\cdots\leq k_1\leq m}Q^{\sum\limits_{j=1}^{p}k_{j-1}k_j-2m\sum\limits_{j=1}^{p}k_j}T^{\sum\limits_{j=1}^{p}k_j}\nonumber\\
            &\quad \times \sqbinom{m}{k_1}_{Q^{-1}}\sqbinom{k_1}{k_2}_{Q^{-1}}\cdots\sqbinom{k_{p-1}}{k_p}_{Q^{-1}}.
        \end{align}
        Note that 
        \begin{align*}
            f(m):
            &=\prod_{k=1}^{m}\frac{Q^{\frac{k}{2}}-Q^{-\frac{k}{2}}}{Q^{\frac{1}{2}}-Q^{-\frac{1}{2}}}\cdot\prod_{k=1}^{m}\frac{1}{1-Q^k}\\
            &=\prod_{k=1}^{m}\left(\frac{Q^{-k/2}(1-Q^k)}{Q^{-1/2}(1-Q)} \cdot \frac{1}{1-Q^k}\right)\\
            &=\frac{Q^{\frac{1}{2}\sum_{k=1}^m (1-k)} }{{(1-Q)}^m}\\
            &=\frac{Q^{\frac{1}{4}(m-m^2)}}{{(1-Q)}^m}.
        \end{align*}
        Hence by Lemma~\ref{cor:recursive-formula-for-bottom-row},
        \begin{align*}
            \dim(\overline{\mathrm{HHH}}_{\wedge^m}(T(2,2p+1))) 
            &= (1-T) {f(m)}^{-1} p'_e(1^m 0^m, 1^m 0^{2pm})\\
            &= Q^{\frac{3m^2}{4}-\frac{3m}{4}}
            \mathop{\sum\sum\cdots\sum}_{0\leq k_p\leq k_{p-1}\leq\cdots\leq k_1\leq m}
            Q^{\sum\limits_{j=1}^{p}k_{j-1}k_j-2m\sum\limits_{j=1}^{p}k_j}T^{\sum\limits_{j=1}^{p}k_j}\\
            &\hspace*{5em}\times\sqbinom{m}{k_1}_{Q^{-1}}\sqbinom{k_1}{k_2}_{Q^{-1}}\cdots\sqbinom{k_{p-1}}{k_p}_{Q^{-1}}
        \end{align*}
        This completes the proof.
    \end{proof}




\printbibliography[title={References}]
\addcontentsline{toc}{section}{\protect\tocsection{}{}{References}}

\end{document}